\documentclass[graybox]{svmult}

\usepackage{mathptmx}       
\usepackage{helvet}         
\usepackage{courier}        
\usepackage{type1cm}        
\usepackage{makeidx}         
\usepackage{graphicx}        
\usepackage{multicol}        
\usepackage[bottom]{footmisc}

\usepackage[utf8]{inputenc}
\usepackage{subfigure}
\usepackage{amsmath,color,amssymb}
\usepackage{psfrag,picture}
\usepackage{tikz}
\usepackage{physics}
\usepackage{bm}
\usepackage{array}
\usepackage{anyfontsize}
\usepackage{t1enc}
\usepackage{color,soul}
\usepackage{multirow}
\usepackage{titlesec}
\newcommand{\vct}[1]{\bm{\mathsf{#1}}}
\newcommand{\pvct}[1]{\bm{#1}}
\newcommand{\mtx}[1]{\bm{\mathsf{#1}}}

\newcommand{\pxx}{\pvct{x}}

\newcommand{\uu}{\vct{u}}
\newcommand{\vv}{\vct{v}}

\newcommand{\TT}{\vct{T}}

\newcommand{\ba}{\begin{array}}
\newcommand{\ea}{\end{array}}
\newcommand{\bd}{\begin{displaymath}}
\newcommand{\ed}{\end{displaymath}}
\newcommand{\bi}{\begin{itemize}}
\newcommand{\ei}{\end{itemize}}
\newcommand{\bn}{\begin{enumerate}}
\newcommand{\en}{\end{enumerate}}

\newcommand{\f}{\frac}

\newcommand{\be}{\begin{enumerate}}
\newcommand{\ee}{\end{enumerate}}

\makeindex     

\begin{document}

\title*{HPS Accelerated Spectral Solvers for Time Dependent Problems}
\titlerunning{HPS Accelerated Solvers for Time Dependent Problems}
\author{Tracy Babb and Per-Gunnar Martinsson and Daniel Appel\"{o}}
\institute{Tracy Babb \at University of Colorado, Boulder, USA, \email{tracy.babb@colorado.edu}
\and Per-Gunnar Martinsson \at University of Texas, Austin, USA, \email{pgm@ices.utexas.edu}
\and Daniel Appel\"{o} \at University of Colorado, Boulder, USA, \email{daniel.appelo@colorado.edu}}
%
%
\maketitle

\abstract{A high-order convergent numerical method for 
solving linear and non-linear parabolic PDEs is presented. The time-stepping is done via an
explicit, singly diagonally implicit Runge-Kutta (ESDIRK) method
of order 4 or 5, and for the implicit solve, we use the recently developed 
``Hierarchial Poincar\'{e}-Steklov (HPS)'' method. The HPS method combines
a multidomain spectral collocation discretization technique (a ``patching method'')
with a nested-dissection type direct solver. In the context under consideration, the
elliptic solve required in each time-step involves the same coefficient matrix, which
makes the use of a direct solver particularly effective. The manuscript describes the methodology and presents numerical experiments.}

\section{Introduction}
\label{sec:intro}
This manuscript describes a highly computationally efficient solver for
equations of the form
\begin{equation}
\kappa \frac{\partial u}{\partial t} = \mathcal{L} u (\pxx, t) + h(u,\pxx,t), \ \  \pxx \in \Omega, t > 0,
\label{eq:parabolic}
\end{equation}
with initial data $u(\pxx,0) = u_{0}(\pxx)$. Here $\mathcal{L}$ is an elliptic operator acting on a fixed domain $\Omega$ and $h$ is lower order, possibly nonlinear terms. We take $\kappa$ to be real or imaginary, allowing for parabolic and Schr\"{o}dinger type equations.
We desire the benefits that can be gained from an implicit solver, such as
L-stability and stiff accuracy, which means that the computational bottleneck
will be the solution of a sequence of elliptic equations set on $\Omega$.
In situations where the elliptic equation to be solved is the same in each
time-step, it is highly advantageous to use a \textit{direct} (as opposed to
\textit{iterative}) solver. In a direct solver, an approximate solution operator
to the elliptic equation is built once. The cost to build it is typically higher
than the cost required for a single elliptic solve using an iterative method such
as multigrid, but the upside is that after it has been built, each subsequent
solve is very fast. In this manuscript, we argue that a particularly efficient
direct solver to use in this context is a method obtained by combining a multidomain
spectral collocation discretization (a so-called ``patching method'', see e.g. Ch. 5.13 in \cite{canuto2007spectral})
with a nested dissection type solver. It has recently been demonstrated
\cite{2016_martinsson_fast_poisson,
2013_martinsson_DtN_linearcomplexity,
2012_spectralcomposite}
that this combined scheme, which we refer to as a ``Hierarchial Poincar\'{e}-Steklov (HPS)'' solver,
can be used with very high local discretization orders (up to $p=20$ or higher) without
jeopardizing either speed or stability, as compared to lower order methods.

In this manuscript, we investigate the stability and accuracy that is obtained
when combining high-order time-stepping schemes with the HPS method for solving elliptic
equations. We restrict attention to relatively simple geometries (mostly rectangles). The
method can without substantial difficulty be generalized to domains that can naturally be
expressed as a union of rectangles, possibly mapped via curvilinear smooth parameter maps. 


\section{The Hierarchical Poincar\'e-Steklov Method}
In this section, we describe a computationally efficient and highly
accurate technique for solving an elliptic PDE of the form
\begin{equation}
\label{eq:basic}
\begin{aligned}
\mbox{}[Au](\pxx) =&\ g(\pxx),\quad &\pxx \in \Omega,\\
   u(\pxx) =&\ f(\pxx),\quad &\pxx \in \Gamma,
\end{aligned} 
\end{equation}
where $\Omega$ is a domain with boundary $\Gamma$, and
where $A$ is a variable coefficient elliptic differential operator
\begin{multline*}
[Au](\pxx) = -c_{11}(\pxx)[\partial_{1}^{2}u](\pxx)
-2c_{12}(\pxx)[\partial_{1}\partial_{2}u](\pxx)
-c_{22}(\pxx)[\partial_{2}^{2}u](\pxx)\\
+c_{1}(\pxx)[\partial_{1}u](\pxx)
+c_{2}(\pxx)[\partial_{2}u](\pxx)
+c(\pxx)\,u(\pxx)
\end{multline*}
with smooth coefficients. In the present context, (\ref{eq:basic}) represents
an elliptic solve that is required in an implicit time-descretization technique
of a parabolic PDE, as discussed in Section \ref{sec:intro}. For simplicity,
let us temporarily suppose that the domain $\Omega$ is rectangular; the extension
to more general domains is discussed in Remark \ref{remark:generaldomain}.

Our ambition here is merely to provide a high level description of the method;
for implementation details, we refer to
\cite{2016_martinsson_fast_poisson,
2013_martinsson_ItI,
2013_martinsson_DtN_linearcomplexity,
2016_martinsson_HPS_3D,
2012_spectralcomposite,
2015_HPS_tutorial}.

\subsection{Discretization} \label{sec:discretization}
\begin{figure}
\begin{center}
\includegraphics[width=60mm]{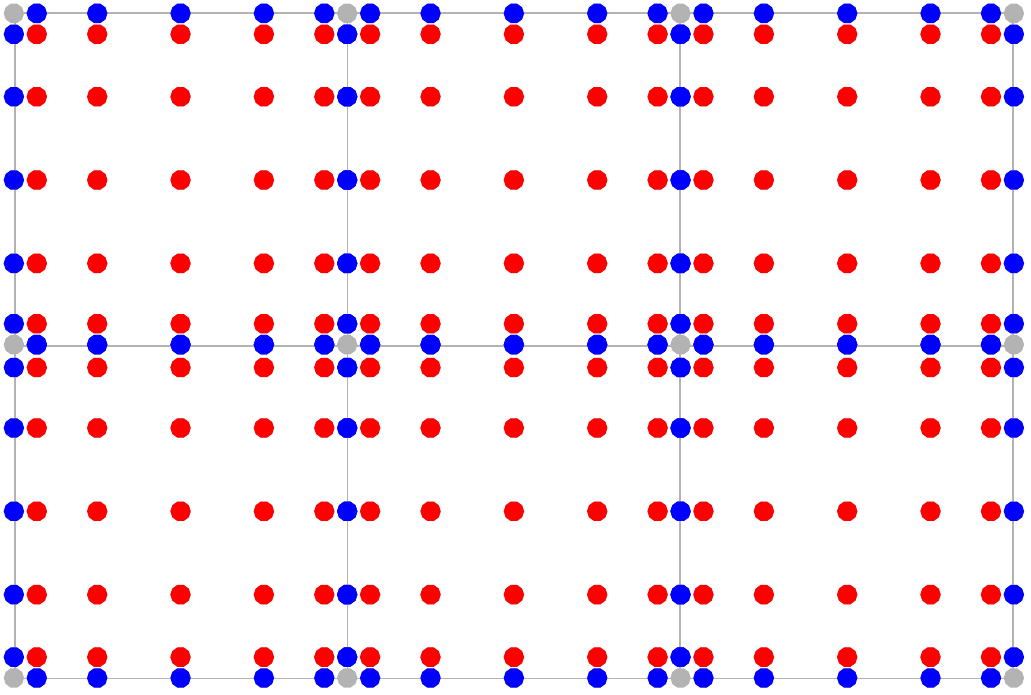}
\end{center}
\caption{The domain $\Omega$ is split into $n_{1} \times n_{2}$ squares, each of size $h\times h$. In the
figure, $n_{1}=3$ and $n_{2} = 2$. Then on each box, a $p\times p$ tensor product grid of Chebyshev nodes
is placed, shown for $p=7$. At red nodes, the PDE (\ref{eq:basic}) is enforced via collocation of the
spectral differentiation matrix. At the blue nodes, we enforce continuity of the normal fluxes. Observe
that the corner nodes (gray) are excluded from consideration.}
\label{fig:full_nodes}
\end{figure}

We split the domain $\Omega$ into $n_{1} \times n_{2}$ boxes, each of size $h \times h$. Then on each box, we place
a $p\times p$ tensor product grid of Chebyshev nodes, as shown in Figure \ref{fig:full_nodes}. We use collocation to
discretize the PDE (\ref{eq:basic}). With $\{\pxx_{i}\}_{i=1}^{N}$ denoting the collocation points, the
vector $\uu$ that represents our approximation to the solution $u$ of (\ref{eq:basic}) is given simply by
$ \uu(i) \approx u(\pxx_{i}).$ We then discretize (\ref{eq:basic}) as follows:
\begin{enumerate}
\item For each collocation node that is \textit{internal} to a box (red nodes in Figure \ref{fig:full_nodes}), we
enforce (\ref{eq:basic}) by directly collocating the spectral differential operator supported on the box, as
described in, e.g., Trefethen \cite{2000_trefethen_spectral_matlab}.
\item For each collocation node on an \textit{edge} between two boxes (blue nodes in Figure \ref{fig:full_nodes}),
we enforce that the normal fluxes across the edge be continuous. For instance, for a node $\pxx_{i}$ on a vertical
line, we enforce that $\partial u/\partial x_{1}$ is continuous across the edge by equating the values for
$\partial u/\partial x_{1}$ obtained by spectral differentiation of the boxes to the left and to the right of the edge.
For an edge node that lies on the external boundary $\Gamma$, simply evaluate the normal derivative at the node, as obtained by spectral
differentiation in the box that holds the node.
\item All \textit{corner} nodes (gray in Figure \ref{fig:full_nodes}) are dropped from consideration. By a happy
coincidence, it turns out that these values do not contribute to any of the spectral derivatives on the remaining nodes
\cite{2017_gillman_HPS_adaptive}.
\end{enumerate}

Since we exclude the corner nodes from consideration, the total number of nodes in the grid equals
$
N = (p-2)\bigl(p\,n_{1}\,n_{2} + n_{1} + n_{2}\bigr) \approx p^{2}\,n_{1}\,n_{2}.
$
The discretization procedure described then results in an $N\times N$ matrix $\mtx{A}$. For a node $i$, the
value of $\mtx{A}(i,:)\uu$ depends on what type of node $i$ is:
$$
\mtx{A}(i,:)\uu \approx \left\{\begin{array}{l}
[Au](\pxx_{i})\ \mbox{for any interior (red) node,}\\
0\ \mbox{for any edge node (blue) not on $\Gamma$,}\\
\partial u /\partial n\ \mbox{for any edge node (blue) on $\Gamma$.}
\end{array}\right.
$$
This matrix $\mtx{A}$ can be used to solve BVPs with a variety of different boundary conditions,
including Dirichlet, Neumann, Robin, and periodic \cite{2012_spectralcomposite}.

In many situations, a simple uniform mesh of the type shown in Figure \ref{fig:full_nodes}
is not optimal, since the regularity in the solution may vary greatly, due to corner singularities,
localized loads, etc. The HPS method can easily be adapted to handle local refinement. The essential
difficulty that arises is that when boxes of different sizes are joined, the collocation nodes
along the joint boundary will not align. It is demonstrated in \cite{2016_martinsson_fast_poisson,2017_gillman_HPS_adaptive}
that this difficulty can stably and efficiently be handled by incorporating local interpolation operators.


\subsection{A Hierarchical Direct Solver}

A key observation in previous work on the HPS method is that the sparse linear
system that results from the discretization technique described in Section \ref{sec:discretization}
is particularly well suited for direct solvers, such as the well-known multifrontal
solvers that compute an LU-factorization of a sparse matrix. The key is to minimize
fill-in by using a so called nested dissection ordering \cite{george_1973,1989_directbook_duff}.
Such direct solvers are very powerful in a situation where a sequence of linear systems
with the same coefficient matrix needs to be solved, since each solve is very fast once
the coefficient matrix has been factorized. This is precisely the environment under
consideration here. The particular advantage of combining the multidomain spectral
collocation discretization described in Section \ref{sec:discretization} is that the
time required for factorizing the matrix is \textit{independent} of the local discretization
order. As we will see in the numerical experiments, this enables us to attain both very
high accuracy, and very high computational efficiency.

\begin{remark}[General domains]
\label{remark:generaldomain}
For simplicity we restrict attention to rectangular domains in this manuscript.
The extension to domains that can be mapped to a union of rectangles via smooth coordinate
maps is relatively straight-forward, since the method can handle variable coefficient
operators \cite[Sec.~6.4]{2012_spectralcomposite}. Some care must be exercised
since singularities may arise at intersections of parameter maps, which may require
local refinement to maintain high accuracy.
\end{remark}

The direct solver described exactly mimics the classical nested dissection method,
and has the same asymptotic complexity of $O(N^{1.5})$ for the ``build'' (or ``factorization'')
stage, and then $O(N\log N)$ cost for solving a system once the coefficient matrix has been
factorized. Storage requirements are also $O(N\log N)$. A more precise analysis of the complexity
that takes into account the dependence on the order $p$ of the local discretization shows
\cite{2016_martinsson_fast_poisson} that
$T_{\rm build} \sim N\,p^{4} + N^{1.5}$, and
$T_{\rm solve} \sim N\,p^{2} + N\log N$.

\section{Time-stepping Methods}

For high-order time-stepping of (\ref{eq:parabolic}), we use the so called 
Explicit, Singly Diagonally Implicit Runge-Kutta (ESDIRK) methods.
These methods have a Butcher diagram with a constant diagonal $\gamma$ and are of the form
\begin{center}
\begin{tabular}{ c | c c c c c c}
        0   & 0           &             &             &          &           & \\
  $2\gamma$ & $\gamma$    & $\gamma$    &             &          &           & \\
  $c_{3}$   & $a_{3,1}$   & $a_{3,2}$   & $\gamma$    &          &           & \\
  $\vdots$  & $\vdots$    & $\vdots$    & $\ddots$    & $\ddots$ &           & \\
  $c_{s-1}$ & $a_{s-1,1}$ & $a_{s-1,2}$ & $a_{s-1,3}$ & $\cdots$ & $\gamma$  & \\
  1         & $b_{1}$     & $b_{2}$     & $b_{3}$     & $\cdots$ & $b_{s-1}$ & $\gamma$ \\
  \hline
            & $b_{1}$     & $b_{2}$     & $b_{3}$     & $\cdots$ & $b_{s-1}$ & $\gamma$ \\
\end{tabular}
\end{center}
ESDIRK methods offer the advantages of stiff accuracy and L-stability. They are
particularly attractive when used in conjunction with direct solvers since the
elliptic solve required in each stage involves the same coefficient matrix
$(I-h\gamma \mathcal{L})$, where $h$ is the time-step. 

In general we split the right hand side of (\ref{eq:parabolic}) into a stiff part, $F^{[1]}$, that will be treated implicitly using ESDIRK methods, and a part, $F^{[2]}$, that will be treated explicitly (with a Butcher table denoted $\hat c$, $\hat A$, and $\hat b$). Precisely we will use the Additive Runge-Kutta (ARK) methods by Carpenter and Kennedy, \cite{Carpenter_Kennedy_ARK}, of order 3,4 and 5.

We may choose to formulate the Runge--Kutta method in terms of either solving for slopes or solving for stage solutions. We denote these the $k_i$ formulation and the $u_i$ formulation, respectively. When solving for slopes the stage computation is
\begin{eqnarray}
k_{i}^{n} &=& F^{[1]}(t_{n}  + c_{i} \Delta t, u^{n} + \Delta t\sum_{j=1}^{s}a_{ij}k_{j}^{n} + \Delta t\sum_{j=1}^{s} \hat{a}_{ij} l_{j}^{n}),\ \ i = 1,\ldots,s,\\
l_{i}^{n} &=& F^{[2]}(t_{n}  + c_{i} \Delta t, u^{n} + \Delta t\sum_{j=1}^{s}a_{ij}k_{j}^{n} + \Delta t\sum_{j=1}^{s} \hat{a}_{ij} l_{j}^{n}), \ \ i = 1,\ldots,s. \label{eq:explicit_slope}
\end{eqnarray}

Note that the explicit nature of (\ref{eq:explicit_slope}) is encoded in the fact that the elements on the diagonal and above in $\hat{A}$ are zero. Once the slopes have been computed the solution at the next time-step is assembled as
\begin{equation}
u^{n+1} = u^{n} + \Delta t\sum_{j=1}^{s}b_{j}k_{j}^{n} + \Delta t\sum_{j=1}^{s} \hat{b}_{j} l_{j}^{n}. \label{eq:rk_update1}
\end{equation}

If the method is instead formulated in terms of solving for the stage solutions the implicit solves take the form
\begin{equation*}
u_{i}^{n} = u^{n} + \Delta t \sum_{j=1}^{s} \Big( a_{ij}F^{[1]}(t_{n}  + c_{j} \Delta t, u_{j}^{n}) + \hat{a}_{ij} F^{[2]}(t_{n}  + c_{j} \Delta t, u_{j}^{n}) \Big ),
\end{equation*}
and the explicit update for $u^{n+1}$ is given by
\begin{equation*}
u^{n+1} = u^{n} + \Delta t \sum_{j=1}^{s} b_{j}(F^{[1]}(t_{n}  + c_{j} \Delta t, u_{j}^{n}) + F^{[2]}(t_{n}  + c_{j} \Delta t, u_{j}^{n})).
\end{equation*}

The two formulations are algebraically equivalent but offer different advantages. For example, when working with the slopes we do not observe (see experiments presented below) any order reduction due to time-dependent boundary conditions (see e.g. the analysis by Rosales et al. \cite{Rosales}).  On the other hand and as discussed in some detail below, in solving for the slopes the HPS framework requires an additional step to enforce continuity.

We note that it is generally preferred to solve for the slopes when implementing implicit Runge-Kutta methods, particularly when solving very stiff problems where the influence of roundoff (or solver tolerance) errors can be magnified by the Lipschitz constant when solving for the stages directly.

\begin{remark}
The HPS method for elliptic solves was previously used in \cite{2014_haut_hyperbolic}, 
which considered a linear hyperbolic equation
\begin{equation*}
\frac{\partial u}{\partial t} = \mathcal{L} u (\pxx, t), \ \ \pxx \in \Omega, t > 0,
\end{equation*}
where $\mathcal{L}$ is a skew-Hermitian operator. The evolution of the numerical solution can be performed by approximating the  propagator $\exp(\tau \mathcal{L}): L^{2}(\Omega)\rightarrow L^{2}(\Omega)$ via a rational approximation
\begin{equation*}
\exp(\tau \mathcal{L}) \approx \sum_{m = - M}^{M} b_{m}(\tau \mathcal{L} - \alpha_{m})^{-1}.
\end{equation*}
If application of $(\tau \mathcal{L} - \alpha_{m})^{-1}$ to the current solution can be reduced to the solution of an elliptic-type PDE it is straightforward to apply the HPS scheme to each term in the  approximation. A drawback with this approach is that multiple operators must be formed and it is also slightly more convenient to time step non-linear equations using the Runge-Kutta methods we use here.
\end{remark}

There are two modifications to the HPS algorithm that are necessitated by the use of ARK time integrators, we discuss these in the next two subsections.

\subsection{Neumann Data Correction in the Slope Formulation}
In the HPS algorithm the PDE is enforced on interior nodes and continuity of the normal derivative is enforced on the leaf boundary. Now, due to the structure of the update formula (\ref{eq:rk_update1}), if at some time $u^n$ has an error component in the null space of the operator that is used to solve for a slope $k_i$, then this will remain throughout the solution process. Although this does not affect the stability of the method it may result in loss of relative accuracy as the solution evolves. As a concrete example consider the heat equation
\begin{equation}
u_{t} =  u_{xx}, \ \  x \in [0,2], \\ t > 0, \label{eq:WBC}
\end{equation}
with the initial data $u(x,0) =  1-|x-1| $, and with homogenous Dirichlet boundary conditions. We discretize this on two leaves which we denote by $\alpha$ and $\beta$.

Now in the $k_{i}$ formulation, we solve several PDEs for the $k_{i}$ values and update the solution as
\begin{equation*}
u^{n+1} = u^{n} + \Delta t  \sum_{j=1}^{s}b_{j}k_{j}^{n}.
\end{equation*}
Here, even though the individual slopes have continuous derivatives the kink in $u^{n}$ will be propagated to $u^{n+1}$. In this particular example we would end up with the incorrect steady state solution $u(x,t) =  1-|x-1| $.

Fortunately, this can easily be mitigated by adding a consistent penalization of the jump in the derivative of the solution during the merging of two leaves (for details see Section 4 in \cite{2016_martinsson_fast_poisson}). That is, if we denote the jump by $[[\cdot]]$ we replace the condition $0 = [[ Tk + h^{k} ]]$ where $Tk$ is the derivative from the homogenous part and $h^{k}$ is the derivative for the particular solution (of the slope) by the condition \mbox{$[[ Tk + h^{k} -  {\Delta t}^{-1} h^{u} ]] = 0$}.
In comparison to \cite{2016_martinsson_fast_poisson} we get the slightly modified merge formula
\begin{equation*}
\mtx{k}_{i,3} =
\bigl(\TT^{\alpha}_{3,3} - \TT^{\beta}_{3,3}\bigr)^{-1}
\bigl( \TT^{\beta }_{3,2}\mtx{k}_{i,2}
      -\TT^{\alpha}_{3,1}\mtx{k}_{i,1}
      +\vct{h}_{3}^{k,\beta}
      -\vct{h}_{3}^{k,\alpha} - \frac{1}{\Delta t}(\vct{h}_{3}^{u,\alpha}
      -\vct{h}_{3}^{u,\beta})
      \bigr),
\end{equation*}
along with the modified equation for the fluxes of the particular solution on the parent box \begin{align*}
\left[\begin{array}{c} \vv_{1} \\ \vv_{2} \end{array}\right] =&\
\left(\left[\begin{array}{ccc}
\TT_{1,1}^{\alpha} & \mtx{0} \\
\mtx{0} & \TT_{2,2}^{\beta }
\end{array}\right] +
\left[\begin{array}{c}
\TT_{1,3}^{\alpha} \\
\TT_{2,3}^{\beta}
\end{array}\right]\,
\bigl(\TT^{\alpha}_{3,3} - \TT^{\beta}_{3,3}\bigr)^{-1}
\bigl[-\TT^{\alpha}_{3,1}\ \big|\ \TT^{\beta}_{3,2}]
\right)
\left[\begin{array}{c} \mtx{k}_{i,1} \\ \mtx{k}_{i,2} \end{array}\right] +\\
&\
\left[\begin{array}{c} \vct{h}^{k,\alpha}_{1} \\ \vct{h}^{k,\beta}_{2} \end{array}\right] +
\left[\begin{array}{c}
\TT_{1,3}^{\alpha} \\
\TT_{2,3}^{\beta}
\end{array}\right]\,
\bigl(\TT^{\alpha}_{3,3} - \TT^{\beta}_{3,3}\bigr)^{-1}
\bigl(\vct{h}_{3}^{\beta}-\vct{h}_{3}^{\alpha} - \frac{1}{\Delta t}(\vct{h}_{3}^{u,\alpha}
      -\vct{h}_{3}^{u,\beta}) \bigr).
\end{align*}

Due to space we must refer to \cite{2016_martinsson_fast_poisson} for a detailed discussion of these equations.
Briefly, $h^{k,\alpha}$ and $h^{k,\beta}$ above denote the spectral derivative on each child's boundary for the particular solution to the PDE for $k_{i}$ and are already present in \cite{2016_martinsson_fast_poisson}. However, $h^{u,\alpha}$ and $h^{u,\beta}$, which denote the spectral derivative of $u^{n}$ on the boundary from each child box, are new additions. 

The above initial data is of course extreme but we note that the problem persists for any non-polynomial initial data with the size of the (stationary) error depending on resolution of the simulation. We further note that the described penalization removes this problem without affecting the accuracy or performance of the overall algorithm.

\begin{remark}
Although for linear constant coefficient PDE it may be possible to project the initial data in a way so that interior affine functions do not cause the difficulty above, for greater generality, we have chosen to enforce the extra penalization throughout the time evolution.
\end{remark}

\begin{remark}
When utilizing the $u_{i}$ formulation in a purely implicit problem we do not encounter the difficulty described above. This is because we enforce continuity of the derivative in $u_{s}^{n}$ when solving
\begin{equation*}
(I - \Delta t \gamma \mathcal{L}) u_{s}^{n}  =  u^{n} + \Delta t \mathcal{L} \Big (\sum_{j=1}^{s-1} a_{sj} u_{j}^{n} \Big) + \Delta t \sum_{j=1}^{s-1} a_{sj} g(x,t_{n}+c_{j} \Delta t),
\end{equation*}
followed by the update $u^{n+1} = u_{s}^{n}$.
\end{remark}

\subsection{Enforcing Continuity in the Explicit Stage}
The second modification is to the first explicit stage in the $k_{i}$ formulation. Solving a problem with no forcing this stage  is simply
\begin{equation*}
k_{1}^{n} = \mathcal{L} (u_{n}).
\end{equation*}

When, for example, $\mathcal{L}$ is the Laplacian, we must evaluate it on all nodes on the interior of the physical domain. This includes the nodes on the boundary between two leafs where the spectral approximation to the Laplacian can be different if we use values from different leaves. The seemingly obvious choice, replacing the Laplacian on the leaf boundary by the average, leads to instability. However, stability can be restored if we enforce $k_{1}^{n} = \mathcal{L} (u_{n})$ on the interior of each leaf and continuity of the derivative across each leaf boundary. Algorithmically, this is straightforward as these are the same conditions that are enforced in the regular HPS algorithm, except in this case we simply have an identity equation for $k_{1}$ on the interior nodes instead of a full PDE.

Although it is convenient to enforce continuity of the derivative using the regular HPS algorithm it can be done in a more efficient fashion by forming a separate system of equations involving only data on the leaf boundary nodes. In a single dimension on a discretization with $n$ leafs this reduces the work associated with enforcing continuity of the derivative across leaf boundary nodes from solving $n \times (p-1) -1$ equations for $n \times (p-1) -1$ unknowns to solving a tridiagonal system of equations  $n-1$ equations for $n-1$ unknowns.


In two dimensions the system is slightly different, but if we have $n \times n$ leafs with $p \times p$ Chebyshev nodes on each leaf then eliminating the explicit equations for the interior nodes reduces the system to $(p-2) \times 2n$ independent tridiagonal systems of $n-1$ equations with $n-1$ unknowns for a total of $(p-2) \times 2n \times (n-1)$ equations with $(p-2) \times 2n \times (n-1)$ unknowns. 


When the $u_{i}$ formulation is used for a fully implicit problem the intermediate stage values still requires us to evaluate $\mathcal{L} u^{n}$, but this quantity only enters through the body load in the intermediate stage PDEs. The explicit first stage in this formulation is simply $u_{1}^{n} = u^{n}$.  Furthermore, while we must calculate
\begin{equation*}
u^{n+1}  =  u^{n} + \Delta t \mathcal{L} \Big (\sum_{j=1}^{s} a_{sj} u_{j}^{n} \Big),
\end{equation*}
this is equivalent to $u_{s}^{n}$ since $b_{j} = {a}_{sj}$ and we simply take $u^{n+1} = u_{s}^{n}$.

When both explicit and implicit terms are present, we proceed differently. Now, the values of $u_{i}^{n}$ look almost identical to the implicit case and we still avoid the problem of an explicit ``solve'' in $u_{1}^{n}$, but we also have
\begin{equation*}
u^{n+1} = u^{n} + \Delta t \sum_{j=1}^{s} b_{j}(F^{[1]}(t_{n}  + c_{j} \Delta t, u_{j}^{n}) + F^{[2]}(t_{n}  + c_{j} \Delta t, u_{j}^{n}))
\end{equation*}
The ESDIRK method has the property that $b_{j} = {a}_{sj}$, but for the explicit Runge-Kutta method we have $b_{j} \neq \hat{a}_{sj}$.  When the explicit operator $F^{[2]}$ does not contain partial derivatives we need not enforce continuity of the derivative and can simply reformulate the method as
\begin{equation*}
u^{n+1} = u_{s}^{n} + \Delta t \sum_{j=1}^{s}(a_{sj} - \hat{a}_{sj})F^{[2]}(t_{n}  + c_{j} \Delta t, u_{j}^{n})
\end{equation*}

\section{Boundary Conditions}
The above description for Runge-Kutta methods does not address how to impose boundary conditions for a system of ODEs resulting from a discretization of a PDE.  In particular, the different formulations incorporate boundary conditions in slightly different ways.

In this work we consider Dirichlet, Neumann, and periodic boundary conditions. For periodic boundary conditions the intermediate stage boundary conditions are enforced to be periodic for both formulations.  As the $k_{i}$ stage values are approximations to the time derivative of $u$, the imposed Dirichlet boundary conditions for $x \in \Gamma$ are $k_{i}^{n} = u_{t}(x,t_{n} + c_{i} \Delta t)$.  When solving for $u_{i}$ one may attempt to enforce boundary conditions using $u_{i} = u(x,t + c_{i} \Delta t), x \in \Gamma$.  However, as demonstrated in part two of this series and discussed in detail in \cite{Rosales}, this results in order reduction for time dependent boundary conditions. 

In the HPS algorithm, Neumann or Robin boundary conditions are mapped to Dirichlet boundary conditions using the linear Dirichlet to Neumann operator as discussed for example in \cite{2016_martinsson_fast_poisson}.

 \section{Time Dependent Boundary Conditions} \label{sec:td}
 This section discusses time-dependent boundary conditions within the two different Runge-Kutta formulations.  In particular, we investigate the order reduction that has been documented in \cite{Rosales} for implicit Runge-Kutta methods and earlier in \cite{Carpenter_RK_Accuracy} for explicit Runge-Kutta methods. 
 
In this first experiment, introduced in \cite{Rosales}, we solve the heat equation in one dimension
\begin{equation}
u_{t} = u_{xx} + f(t),\qquad  x \in [0,2], \ \ t>0. \label{eq:timeBCEx} 
%
\end{equation}
We set the initial data, Dirichlet boundary conditions and the forcing $f(t)$ so that exact solution is $u(x,t) = \cos(t)$. This example is designed to eliminate the effect of the spatial discretization, with the solution being constant in space and allows for the study of possible order reduction near the boundaries. 

We use the HPS scheme in space and use 32 leafs with $p = 32$ Chebyshev nodes per leaf.  We apply the third, fourth, and fifth order ESDIRK methods from \cite{Carpenter_Kennedy_ARK}. We consider solving for the intermediate solutions, or as we refer to it below ``the $u_{i}$ formulation'' with the boundary condition enforced as \mbox{$u_{i}^{n} = \cos(t_{n} + c_{i}\Delta t)$}.  We also consider solving for the stages, which we refer to as ``the $k_{i}$ formulation'' with boundary conditions imposed as $k_{i}^{n} = -\sin(x,t_{n} + c_{i} \Delta t)$.

\begin{figure}[ht]
\begin{center}
\subfigure[]{\includegraphics[height=.38\textwidth]{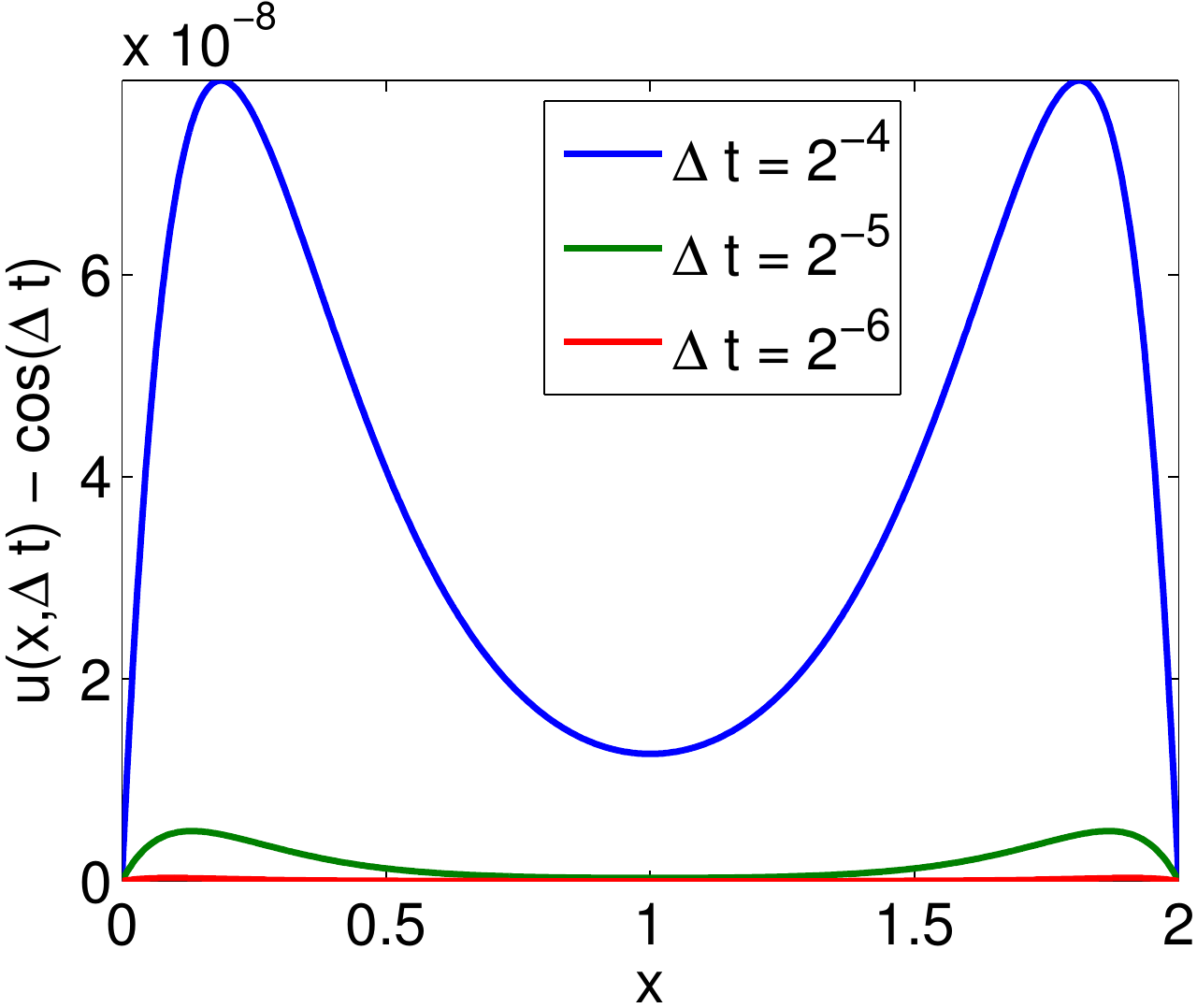}}
\subfigure[] {\includegraphics[height=.38\textwidth]{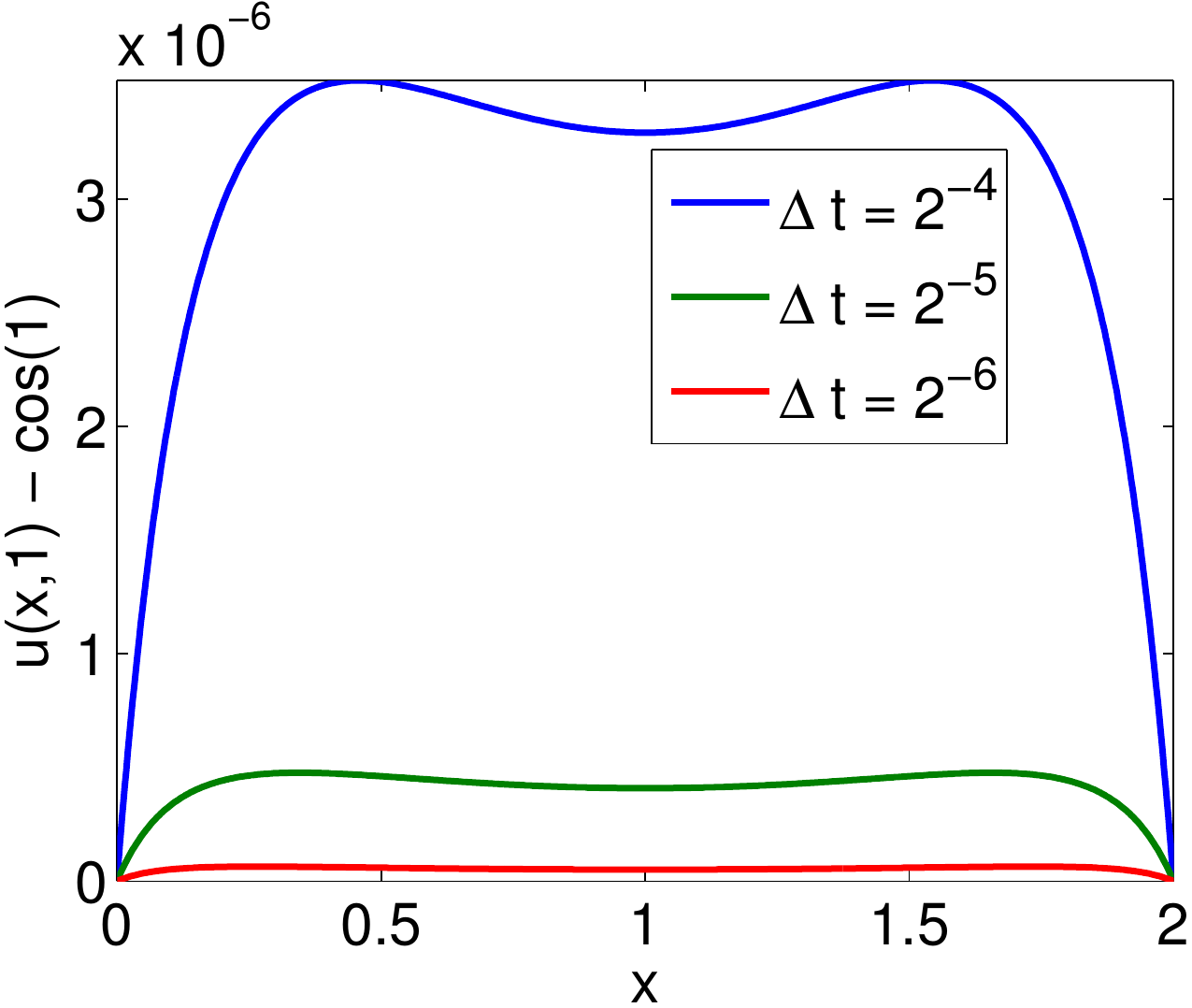}}
\caption{The error in solving (\ref{eq:timeBCEx}). Results are for a 3rd order ESDIRK. Figure (a) displays the single step error which converges with 4th order of accuracy. Figure (b) displays the global error at $t=1$ converging at 3rd order. Both errors converge at one order higher than what is expected from the analysis in \cite{Rosales}.\label{fig:Error_ESDIRKU_O3_x}} 
\end{center}
\end{figure}

Error reduction for time dependent boundary conditions has been studied both in the context of explicit Runge-Kutta methods in e.g. \cite{Carpenter_RK_Accuracy} and more recently for implicit Runge-Kutta methods in \cite{Rosales}. In \cite{Rosales} the authors report observed orders of accuracy equal to two (for the solution $u$) for DIRK methods of order 2, 3, and 4 for the problem (\ref{eq:timeBCEx}) discretized with a finite difference method on a fine grid (the spatial errors are zero) using the $u_i$ formulation.

\begin{figure}[htb]
\begin{center}
\subfigure[]{\includegraphics[height=.38\textwidth]{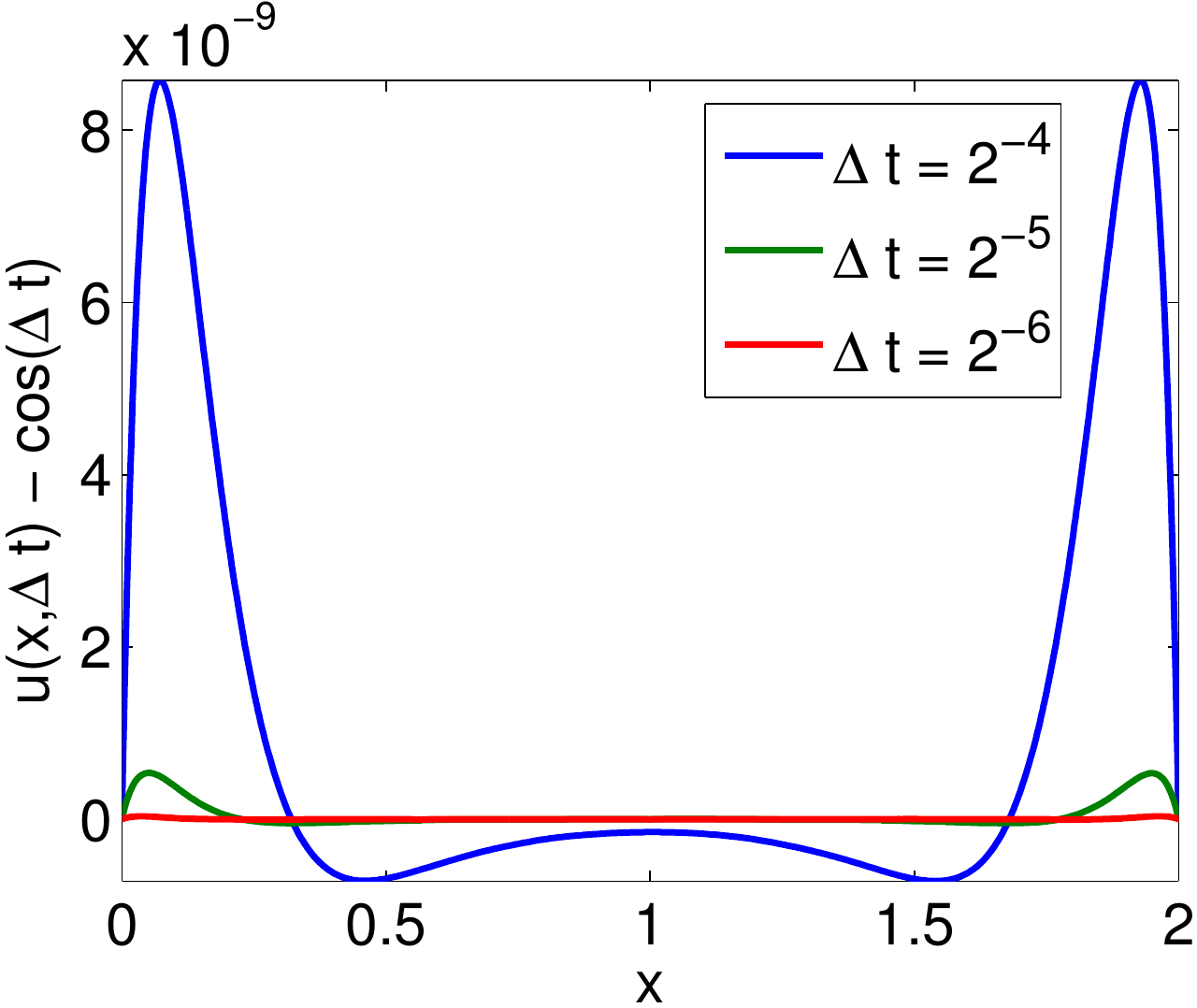}}
\subfigure[] {\includegraphics[height=.38\textwidth]{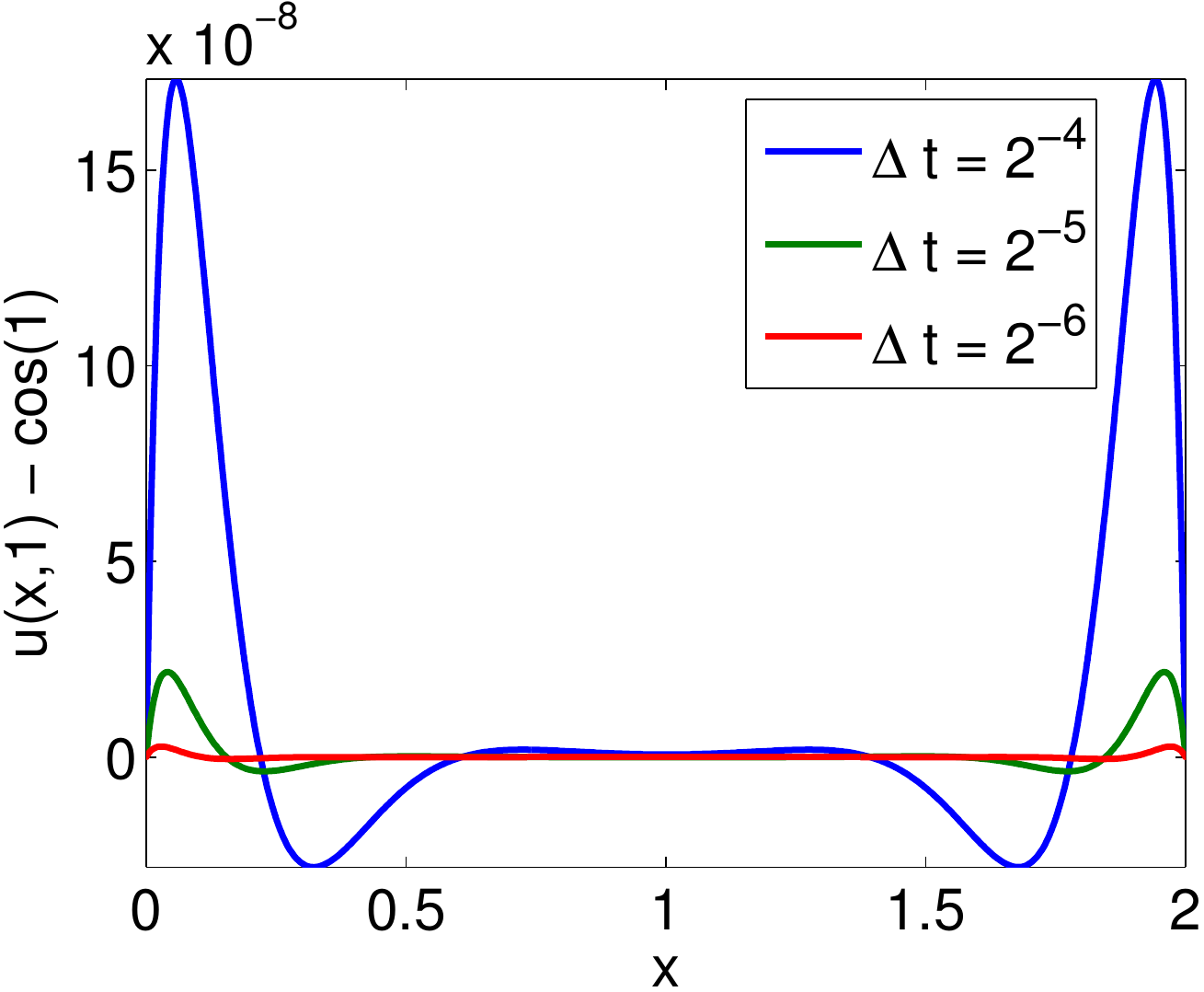}}
\caption{The error in solving (\ref{eq:timeBCEx}). Results are for a 5th order ESDIRK. Figure (a) displays the single step error which converges with 4th order of accuracy. Figure (b) displays the global error at $t=1$ converging at 3rd order. Both errors converge at one order higher than what is expected from the analysis in \cite{Rosales} but still lower than expected. \label{fig:Error_ESDIRKU_O5_x}}
\end{center}
\end{figure}

\begin{figure}[htb]
\begin{center}
\subfigure[]{\includegraphics[height=.383\textwidth]{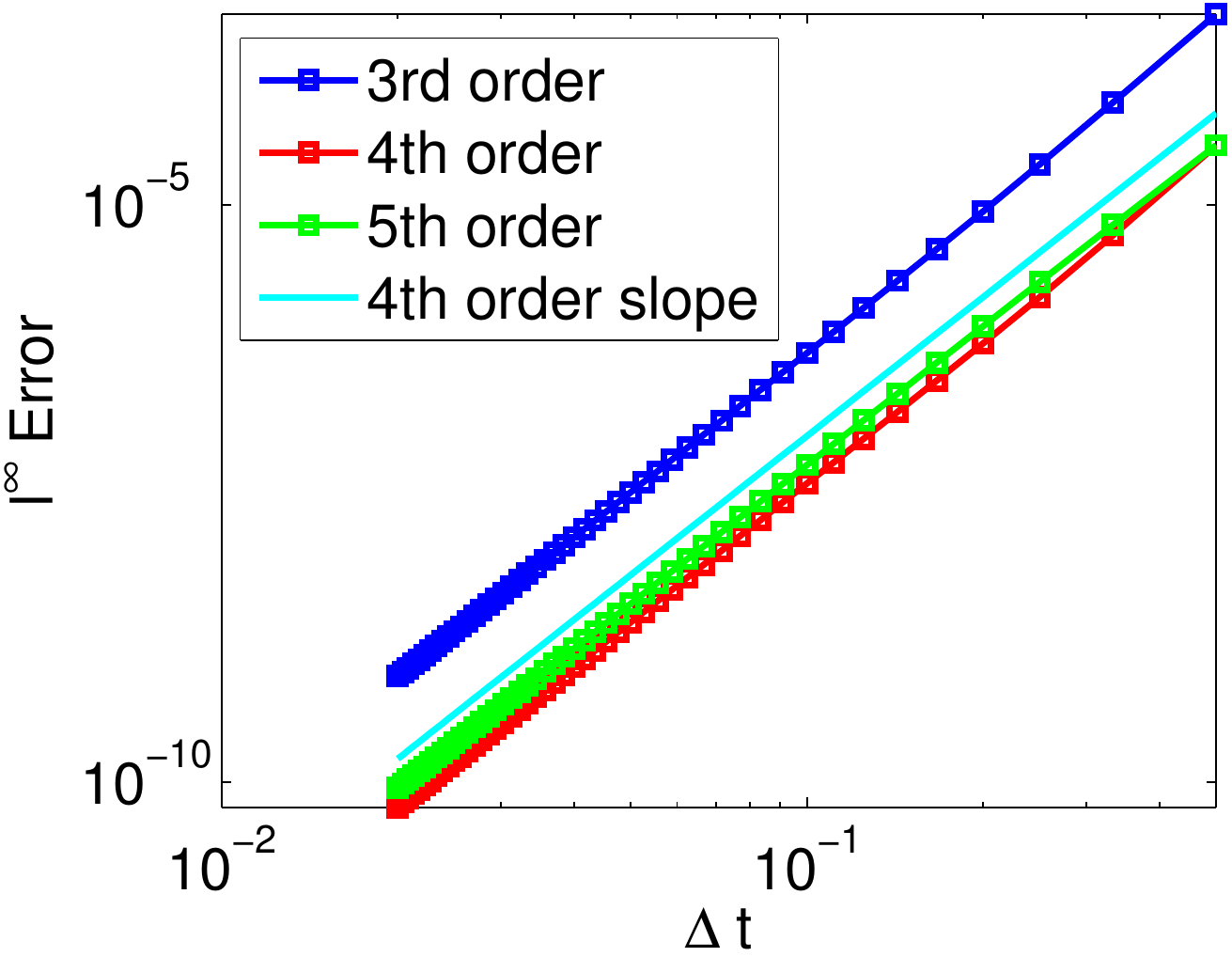}}
\subfigure[] {\includegraphics[height=.383\textwidth]{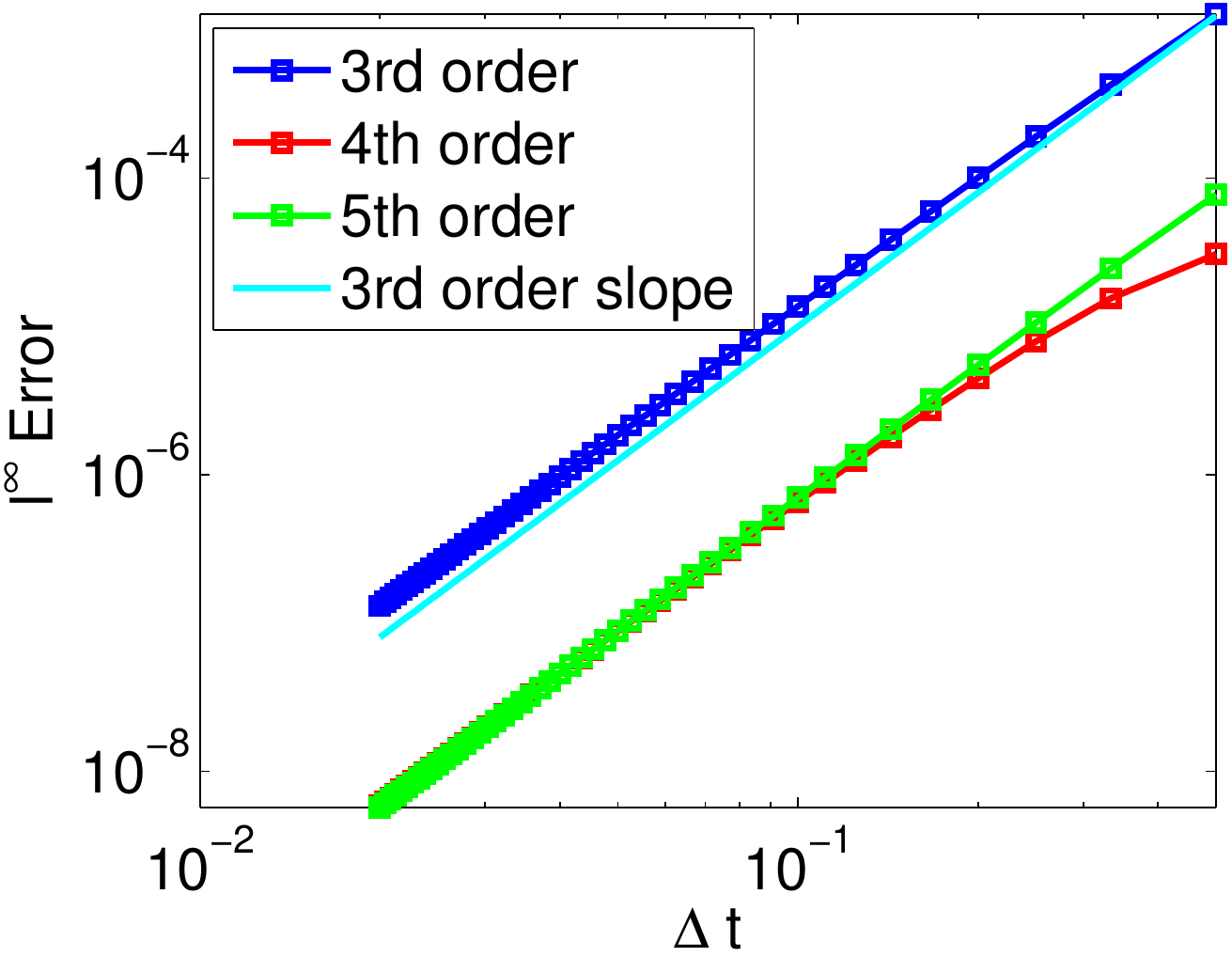}}
\subfigure[]{\includegraphics[height=.35\textwidth]{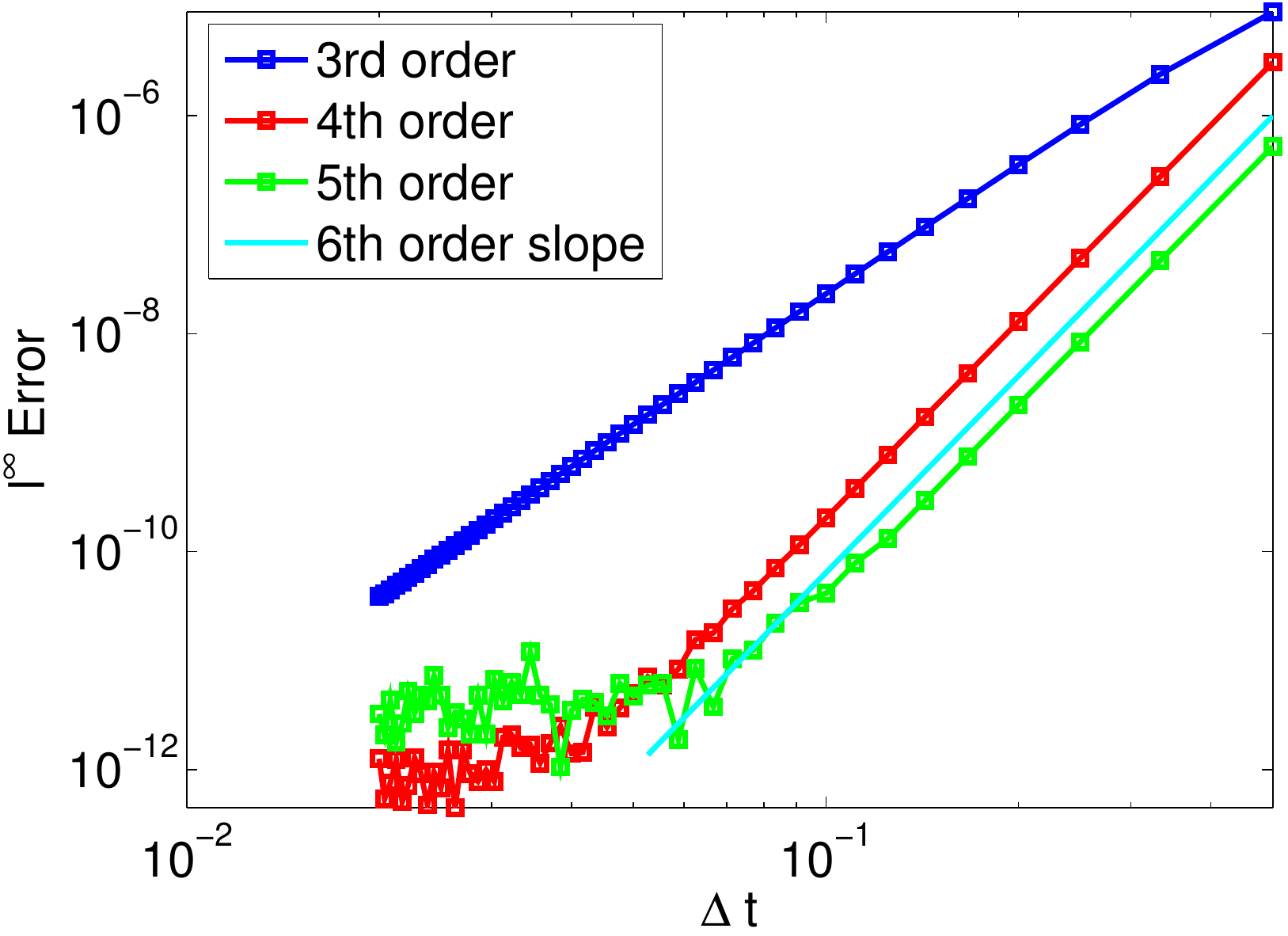}}
\subfigure[] {\includegraphics[height=.35\textwidth]{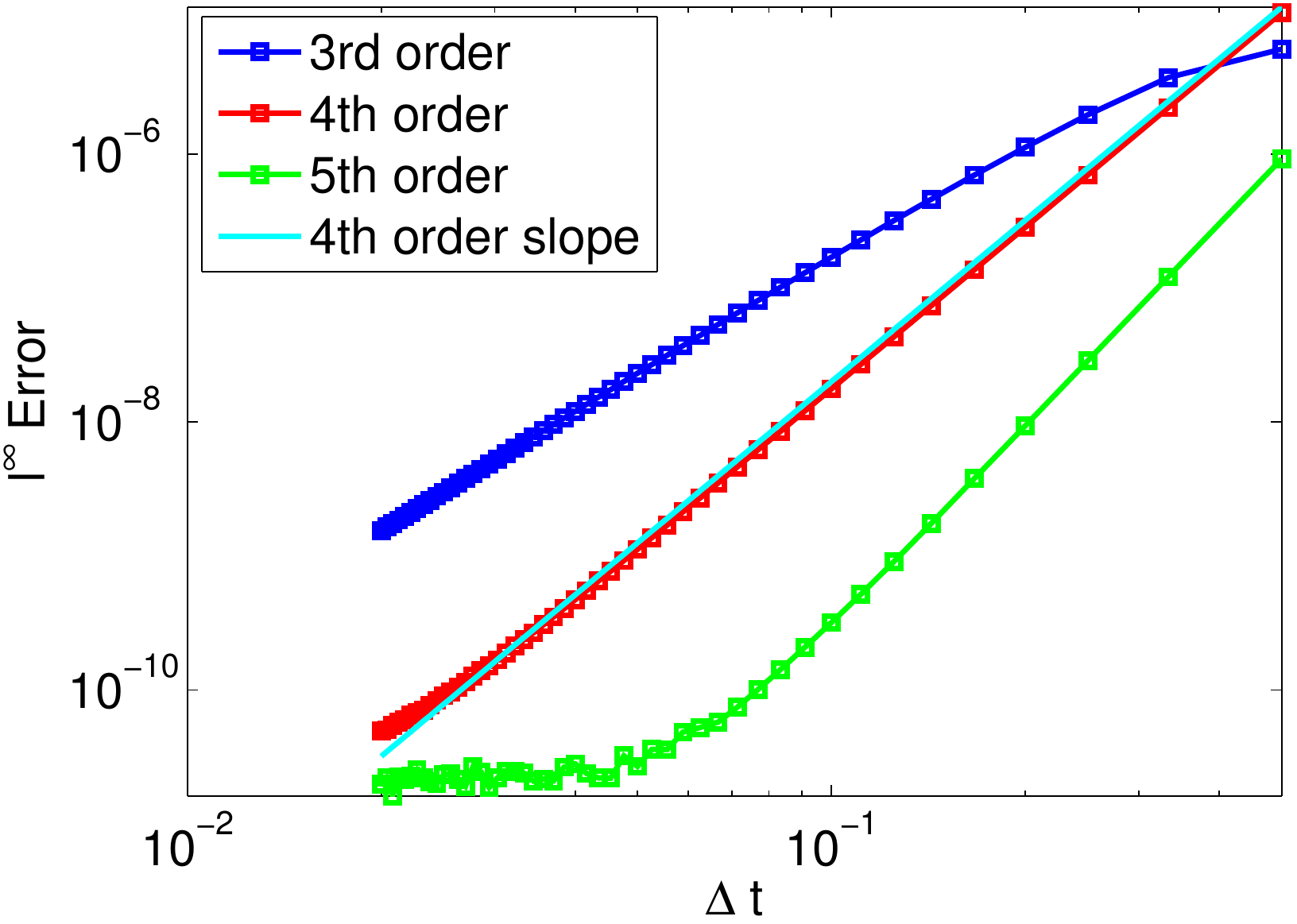}}
\caption{The error in solving (\ref{eq:timeBCEx}) for the 3rd, 4th, and 5th order ESDIRK methods for a sequence of decreasing time steps. Figures (a) and (c) are errors after one time step and (b) and (c) are the errors at time $t=1$. The top row are for the $u_i$ formulation and the bottom row is for the $k_i$ formulation. Note that the $k_i$ formulation is free of order reduction. \label{fig:Error_ESDIRKU_k}}
\end{center}
\end{figure}

Figures \ref{fig:Error_ESDIRKU_O3_x} and \ref{fig:Error_ESDIRKU_O5_x} show the error for the third and fifth order ESDIRK methods, respectively, as a function of $x$ for a single step and at the final time $t=1$. 
Figure \ref{fig:Error_ESDIRKU_k} shows the maximum error for the third, fourth, and fifth order methods as a function of time step $\Delta t$ after a single step and at the final time $t=1$.

In general, for a method of order $p$ we expect that the single step error decreases as $\Delta t^{p+1}$ while the global error decreases as $\Delta t^{p}$.  However, with time dependent boundary conditions implemented as $u_{i}^{n} = \cos(t_{n} + c_{i}\Delta t)$ the results in \cite{Rosales} indicate that the rate of convergence will not exceed two for the single step or global error.  

The results for the third order method ($p=3$) displayed in Figure \ref{fig:Error_ESDIRKU_O3_x} show that the single step error decreases as $\Delta t^{p+1}$ while the global error decreases as $\Delta t^{p}$, which is better than the results documented in \cite{Rosales}. However, we still see that a boundary layer appears to be forming, but it is of the same order as the error away from the boundary. The results for the fifth order method ($p=5$) displayed in Figure \ref{fig:Error_ESDIRKU_O5_x} show that the single step error decreases as $\Delta t^{4}$ while the global error decreases as $\Delta t^{3}$, which is still better than the results documented in \cite{Rosales}.  However, the boundary layer is giving order reduction from $\Delta t^{p+1}$ for the single step error and $\Delta t^{p}$ for the global error. We note that our observations differ from those in \cite{Rosales} but that this possibly can be attributed to the use of a ESDIRK method rather than a DIRK method. 

We repeat the experiment but now we use the $k_{i}$ formulation for Runge-Kutta methods and for the boundary condition we enforce $k_{i}^{n} = -\sin(t_{n} + c_{i}\Delta t)$.  The intuition here is that $k_{i}^{n}$ is an approximation to $u_{t}$ at time $t_{n} + c_{i}\Delta t$ and we use the value of $u_{t}$ for the boundary condition of $k_{i}^{n}$. Intuitively we expect that the fact that we reduce the index of the system of differential algebraic equation in the $u_i$ formulation by differentiating the boundary conditions can restore the design order of accuracy. 

In the previous examples the Runge-Kutta method introduced an error on the interior while the solution on the boundary was exact.  If the error on the boundary is on the same order of magnitude as the error on the interior then the error in $u_{xx}$ is of the correct order, but when the value of $u$ is exact on the boundary it introduces a larger error in $u_{xx}$.  In the $k_i$ formulation, for each intermediate stage we find $u_{xx} = 0$ and then \mbox{$k_{i}^{n} = -\sin(t_{n} + c_{i}\Delta t)$} on the interior and on the boundary.  So at a fixed time the solution is constant in $x$ and a boundary layer does not form. Additionally, the error is constant in $x$ at any fixed time and for a method of order $p$ we obtain the expected behavior where the single step error decreases as $\Delta t^{p+1}$ and the global error decreases as $\Delta t^{p}$. 

Figure \ref{fig:Error_ESDIRKU_k} shows the maximum error for the third, fourth, and fifth order methods as a function of time step $\Delta t$ after a single step and at the final time $t=1$.  The results show that the methods behave exactly as we expect.  The single step error behaves as $\Delta t^{p+1}$ for the 3rd and 5th order methods and $\Delta t^{p+2}$  for the 4th order method.  The 4th order method gives 6th order error in a single step because the exact solution is $u(x,t) = \cos(t)$, which has every other derivative equal to zero at $t=0$ and for a single step we start at $t=0$. The global error behaves as $\Delta t^{p}$ for each method.

\section{Schr\"odinger Equation} \label{sec:sch}
Next we consider the Schr\"{o}dinger equation for $u = u(x,y,t)$ 
\begin{eqnarray}
\begin{aligned}  \label{S2D}
i \hbar  u_t &= - \f{\hbar^2}{2 M} \Delta u+ V(x,y) u, \ \ 
 t>0, \ \ (x,y) \in [x_{ l}, x_{ r}] \times [y_{ b}, y_{ t}],\\
u(x,y,0) &= u_0(x,y).  
\end{aligned}
\end{eqnarray}
Here we nondimensionalize in a way equivalent to setting $M = 1, \hbar = 1$ in the above equation. We choose the potential to be the the harmonic potential
\[
V(x,y) = \f{1}{2} \left( x^2 + y^2 \right).
\] 
This leads to an exact solution
\begin{equation} \label{eq:2DHpot}
u(x,y,t) = A e^{-it} e^{-\f{(x^2+y^2)}{2}}, 
\end{equation}
where we set $A = 1/\sqrt{\sqrt{\pi}}$ and solve until $t=2 \pi$ on the domain $(x,y) \in [-8, 8]^2$. 

\begin{figure}[ht]
 \begin{center}
   \includegraphics[width=.47\textwidth]{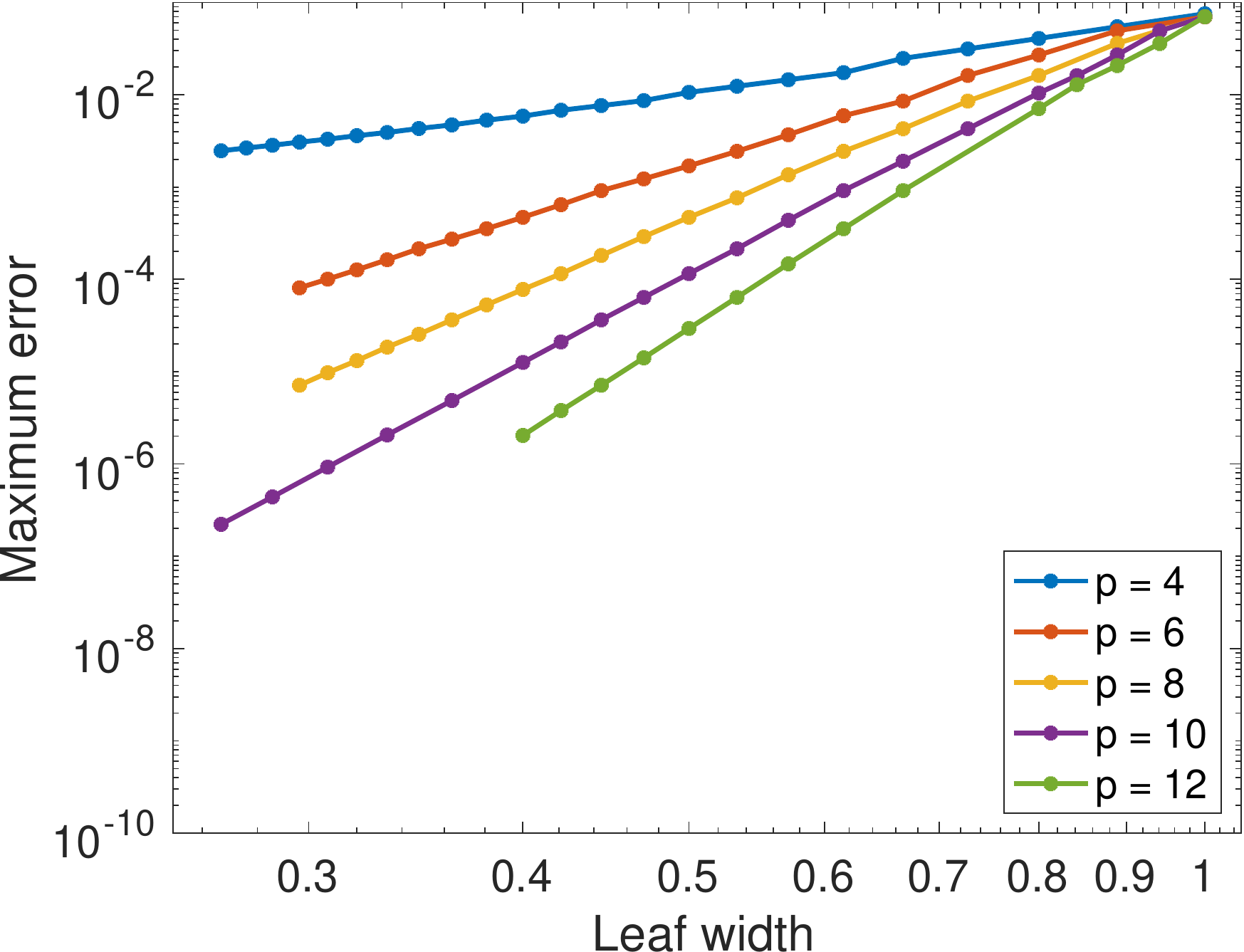}\ \ \ \  \ \ 
   \includegraphics[width=.47\textwidth]{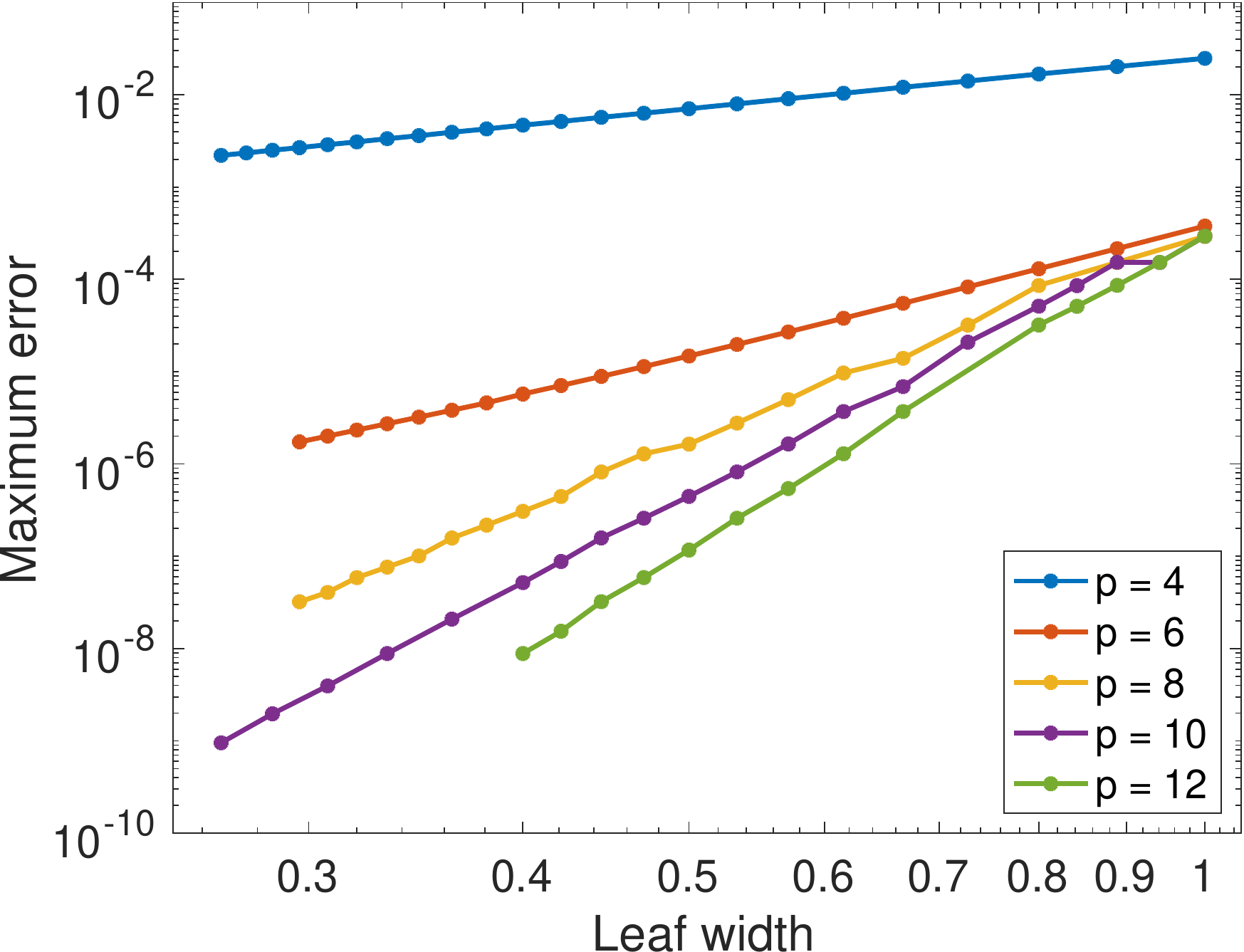}
   \caption{Error in the Schr\"odinger equation as a function of leaf size. The exact solution is (\ref{eq:2DHpot}). \label{fig:Sch_Err}}
\end{center}
\end{figure}
The computational domain is subdivided into $n_{x} \times n_{y}$ panels with $p \times p$ points on each panel.  To begin, we study the order of accuracy with respect to leaf size.  To eliminate the effect of time-stepping errors we scale $\Delta t = h^{p/q_{\rm RK}}$, where $q_{\rm RK}$ is the order of the Runge-Kutta method.  In Figure \ref{fig:Sch_Err} we display the errors as a function of the leaf size for $p = 4, 6, 8, 10, 12, 16$ and for the third and fifth order Runge-Kutta methods ($q_{RK} = 3,5$).  The rates of convergence are found for all three Runge-Kutta methods and summarized in Table \ref{Tab:SCH_rates}. As can be seen from the table, $p=4$ appears to converge at second order, while for higher $p$ we generally observe a rate of convergence approaching to $p$. 
\begin{table}[h]
\begin{center}
\caption{Estimated rates of convergence for different Runge Kutta methods and different orders of approximation. \label{Tab:SCH_rates}}
\begin{tabular}{|l|c|c|c|c|c|c|}
\hline
$p$ & 4    & 6    & 8    & 10    & 12     \\ 
\hline
\hline
ESDIRK3 & 2.59 & 5.73 & 7.72 & 9.69 & 11.47  \\ 
ESDIRK4 & 1.89 & 6.47 & 7.82 & 9.76 & 11.69  \\ 
ESDIRK5 & 1.84 & 4.42 & 7.69 & 9.71 & 11.48  \\ 
\hline
\end{tabular}
\end{center}
\end{table} 

In this problem the efficiency of the method is limited by the order of the Runge-Kutta methods. However, as our methods are unconditionally stable we may enhance the efficiency by using Richardson extrapolation to achieve a highly accurate solution in time.  We solve the same problem, but now we fix $p=12$ and take $5 \cdot 2^n$ time steps, with $n = 0, 1, \dots, 5$.  For the third order ESDIRK method we use $60 \times 60$ leaf boxes.  For the fourth order ESDIRK method we use $90 \times 90$ leaf boxes.  For the fifth order ESDIRK method we use $120 \times 120$ leaf boxes.  Table \ref{Tab:SCH_rates_RE} shows that we can easily achieve much higher accuracy by using Richardson extrapolation.

\begin{table}[ht]
\begin{center}
\caption{Estimated errors at the final time after Richardson extrapolation.  \label{Tab:SCH_rates_RE}}
\begin{tabular}{|l|c|c|c|c|c|c|c|}
\hline
$q_{RK}$ / Extrapolations & 0 & 1 & 2 & 3 & 4 & 5 & 6\\
\hline
3 & 1.32(-1) & 1.01(-2) & 1.27(-4) & 1.17(-5) & 6.98(-8) & 8.62(-10) & 7.40(-6)  \\ 
\hline
4 & 2.70(-4) & 6.46(-6) & 1.23(-7) & 2.95(-10) & 1.59(-11) & 3.70(-14) & 1.20(-11)  \\ 
\hline
5 & 1.28(-3) & 9.67(-6) & 6.30(-8) & 1.86(-10) & 4.11(-13) & 9.27(-14) & 5.08(-11)  \\ 
\hline
\end{tabular}
\end{center}
\end{table}

\begin{table}[ht]
\begin{center}
\caption{Errors computed against a $p$ and $h$ refined solution. The errors are maximum errors at the final time $t=4$. \label{Tab:SCH_rates_EX2}}
\begin{tabular}{|l|c|c|c|c|c|}
\hline
$p$ / Panels & 2 & 4 & 8 & 16 & 32 \\
\hline
8 & 1.11(0) & 1.39(-1) & 8.74(-3) & 1.50(-4) & 2.45(-6)\\ 
\hline
rate & $\ast$ & 3.00 & 3.99 & 5.87 & 5.92 \\ 
\hline
10 & 5.87(-1) & 3.16(-2) & 4.62(-4) & 6.17(-6) & 5.21(-8)\\ 
\hline
rate & $\ast$ & 4.21 & 6.10 & 6.22 & 6.89  \\ 
\hline
\end{tabular}
\end{center}
\end{table} 

Finally, we solve a problem without an analytic solution. In this problem the initial data 
\[
u(x,y,t) = 3 \sin(x) \sin(y) e^{-(x^2+y^2)},
\]
interacts with the weak and slightly non-symmetric potential 
\[
V(x,y) = 1-e^{-(x+0.9y)^4},
\]
allowing the solution to reach the boundary where we impose homogenous Dirichlet conditions. 

We evolve the solution until time $t=4$ using $p = 8$ and 10 and $2, 4, 8, 16$ and 32 leaf boxes in each direction of a domain of size $12 \times 12$. The errors computed against a reference solution with $p=12$ and with 32 leaf boxes can be found in Table \ref{Tab:SCH_rates_EX2}. 
\begin{figure}[htb]
\begin{center}
\subfigure[]{\includegraphics[height=.36\textwidth]{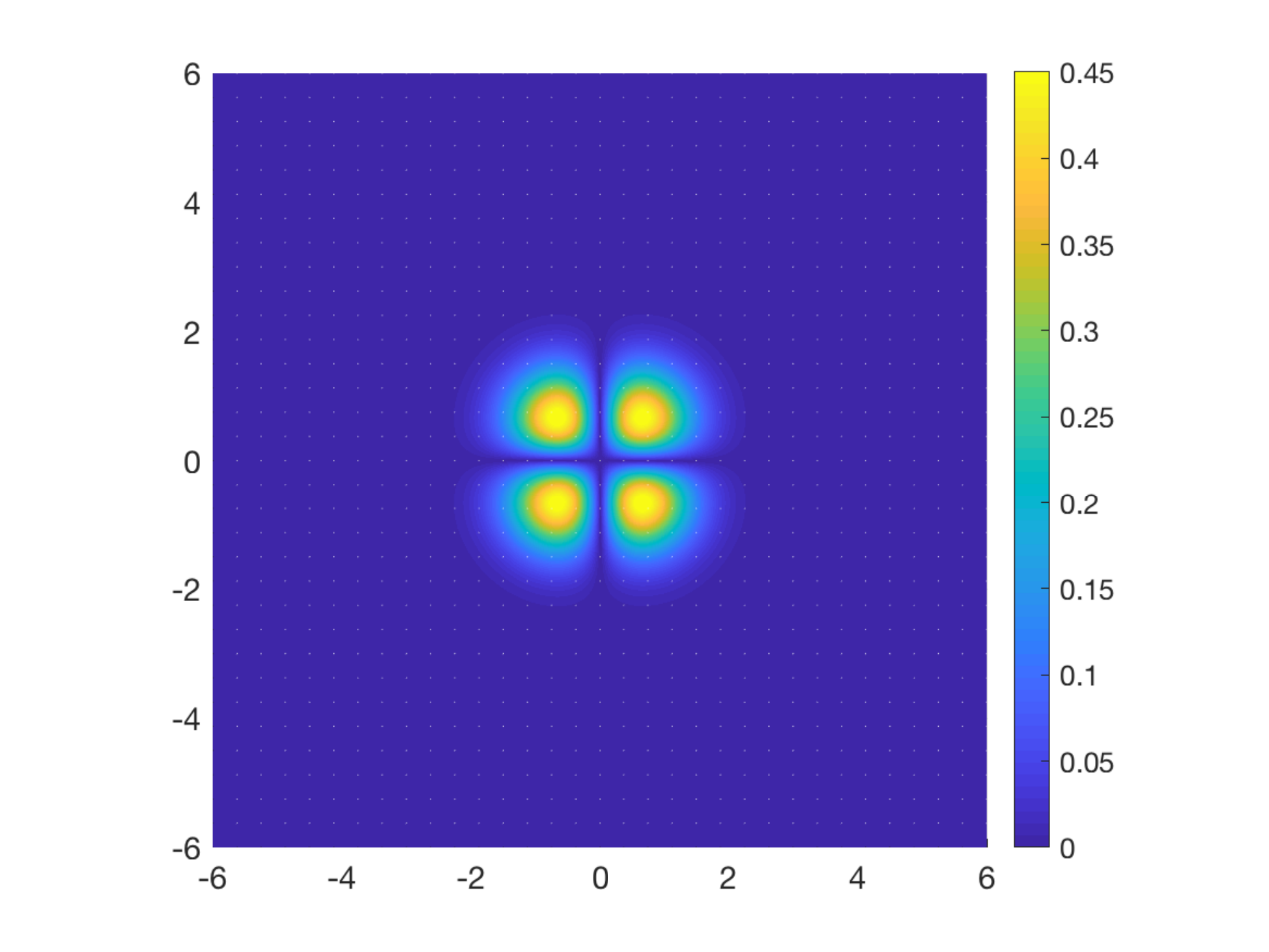}}
\subfigure[] {\includegraphics[height=.36\textwidth]{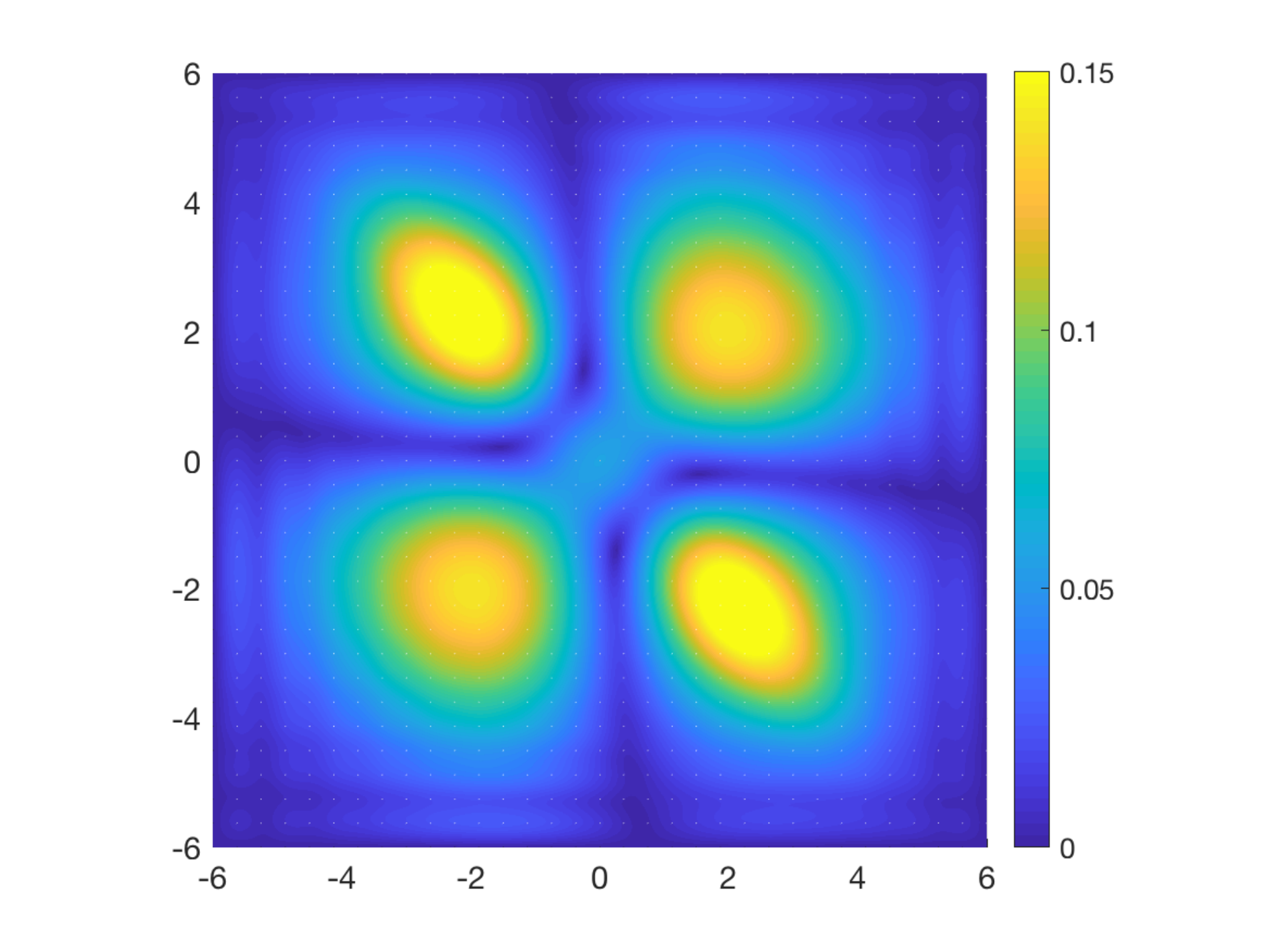}}
\subfigure[]{\includegraphics[height=.36\textwidth]{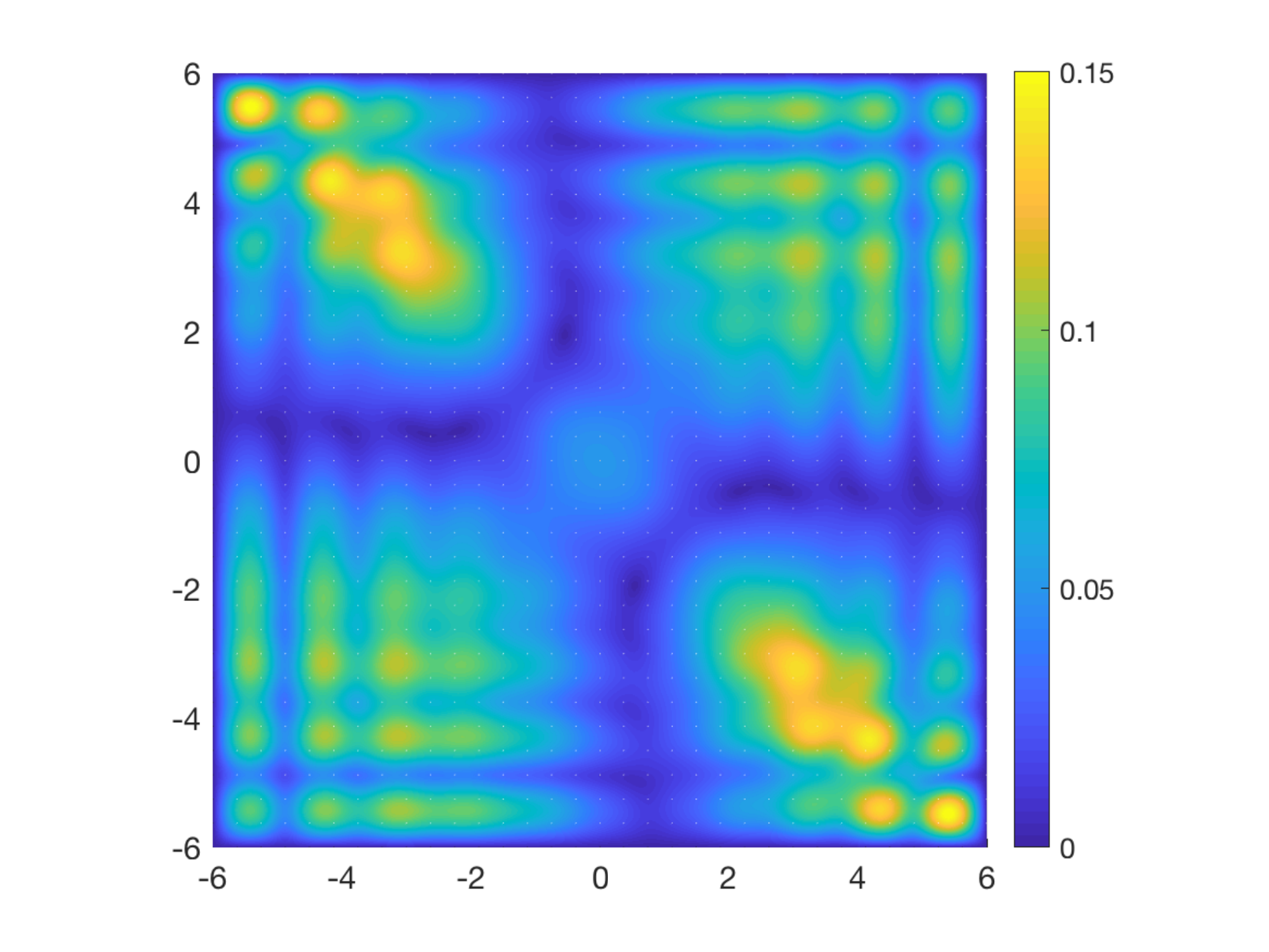}}
\subfigure[] {\includegraphics[height=.36\textwidth]{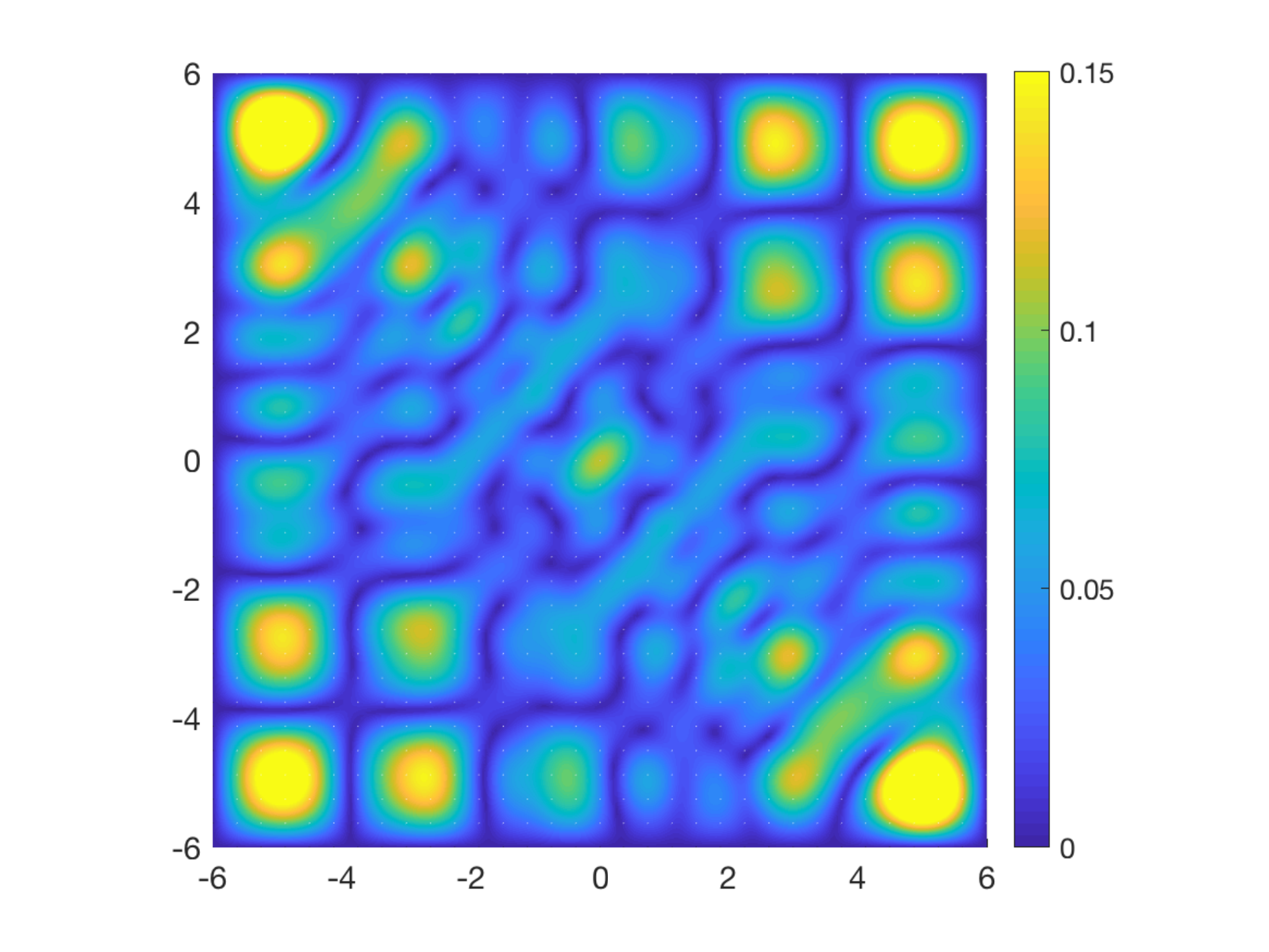}}
\caption{Snapshots of the magnitude of the solution at the initial time $t=0$, the intermediate times $t \approx 1.07$, $t \approx 1.68$ and at the final time $t = 4.0$. \label{fig:SCH_EX2}}
\end{center}
\end{figure}

In Figure \ref{fig:SCH_EX2} we display snapshots of the magnitude of the solution at the initial time $t=0$, the intermediate times $t \approx 1.07$, $t \approx 1.68$ and at the final time $t = 4.0$.



\section{Burgers' Equation in Two Dimensions} \label{sec:bur}
As a first step towards a full blown flow solver we solve Burgers' equation in two dimensions using the additive Runge-Kutta methods described in the first part of this paper. Precisely, we solve the system
\begin{equation}
\mathbf{u}_t + \mathbf{u} \cdot \nabla \mathbf{u} = \varepsilon \Delta \mathbf{u}, \ \ \ \ \mathbf{x} \in [-\pi,\pi]^2, \ \ t>0,
\end{equation}
where $\mathbf{u} = [u(x,y,t),v(x,y,t)]^T$ is the vector containing the velocities in the $x$ and $y$ directions. 

\begin{figure}[!ht]
 \begin{center}
    \includegraphics[width=.3\textwidth]{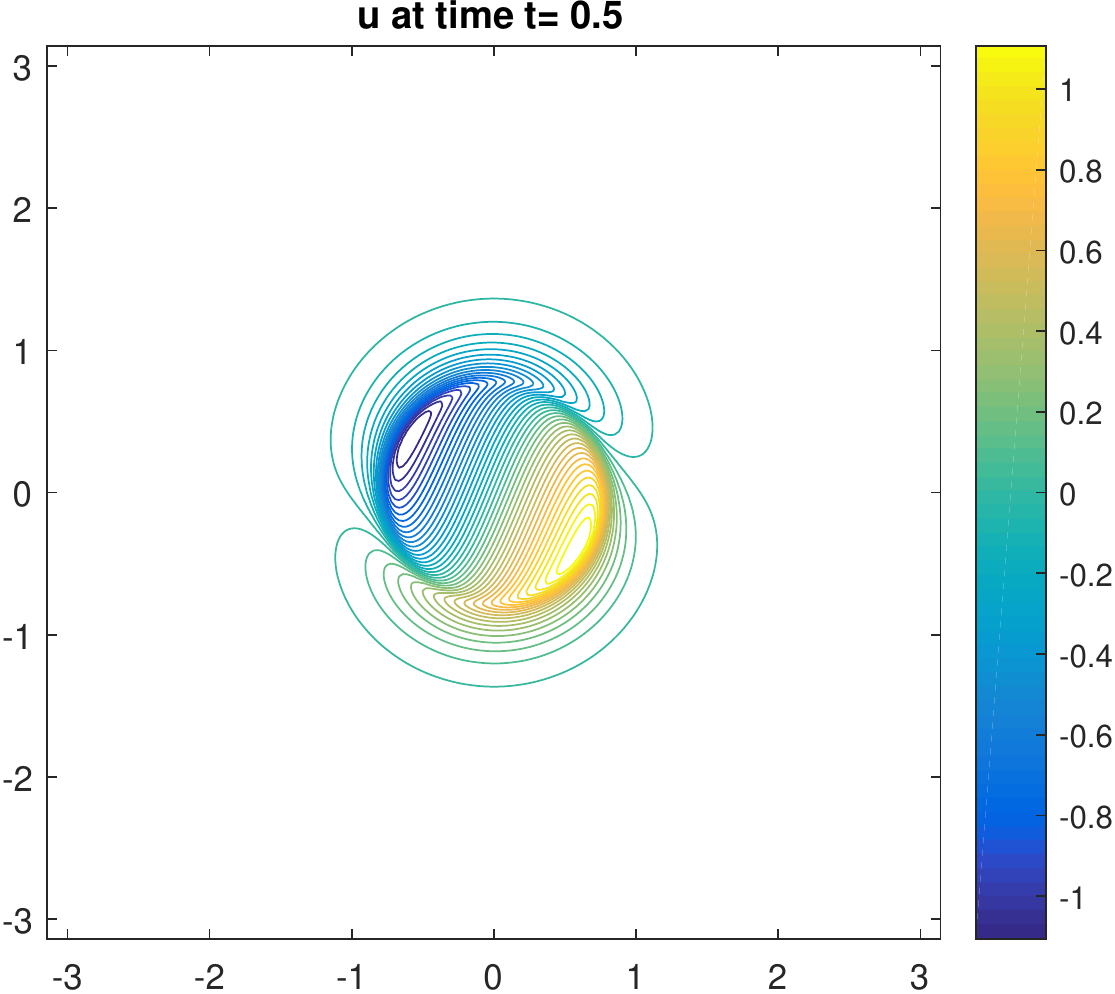}
   \includegraphics[width=.3\textwidth]{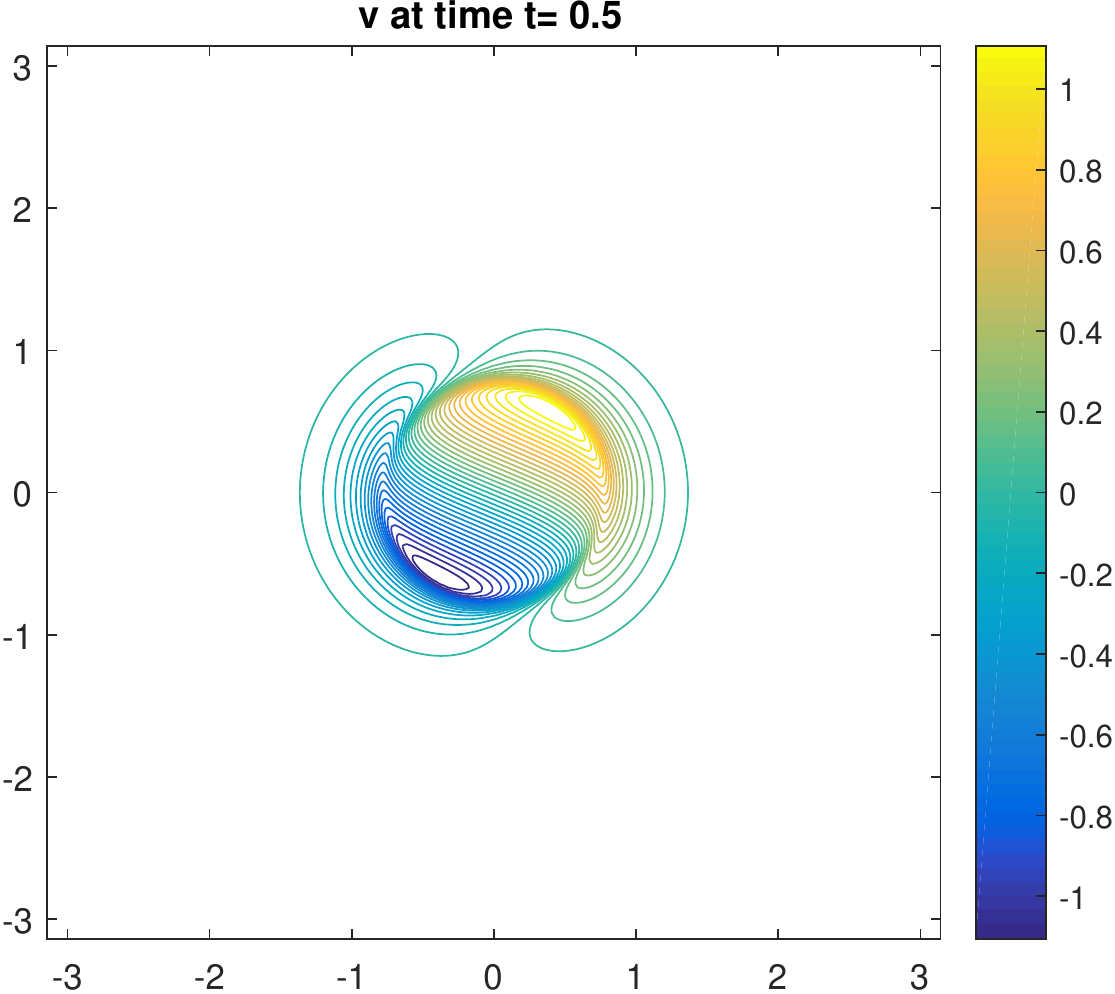}
      \includegraphics[width=.335\textwidth]{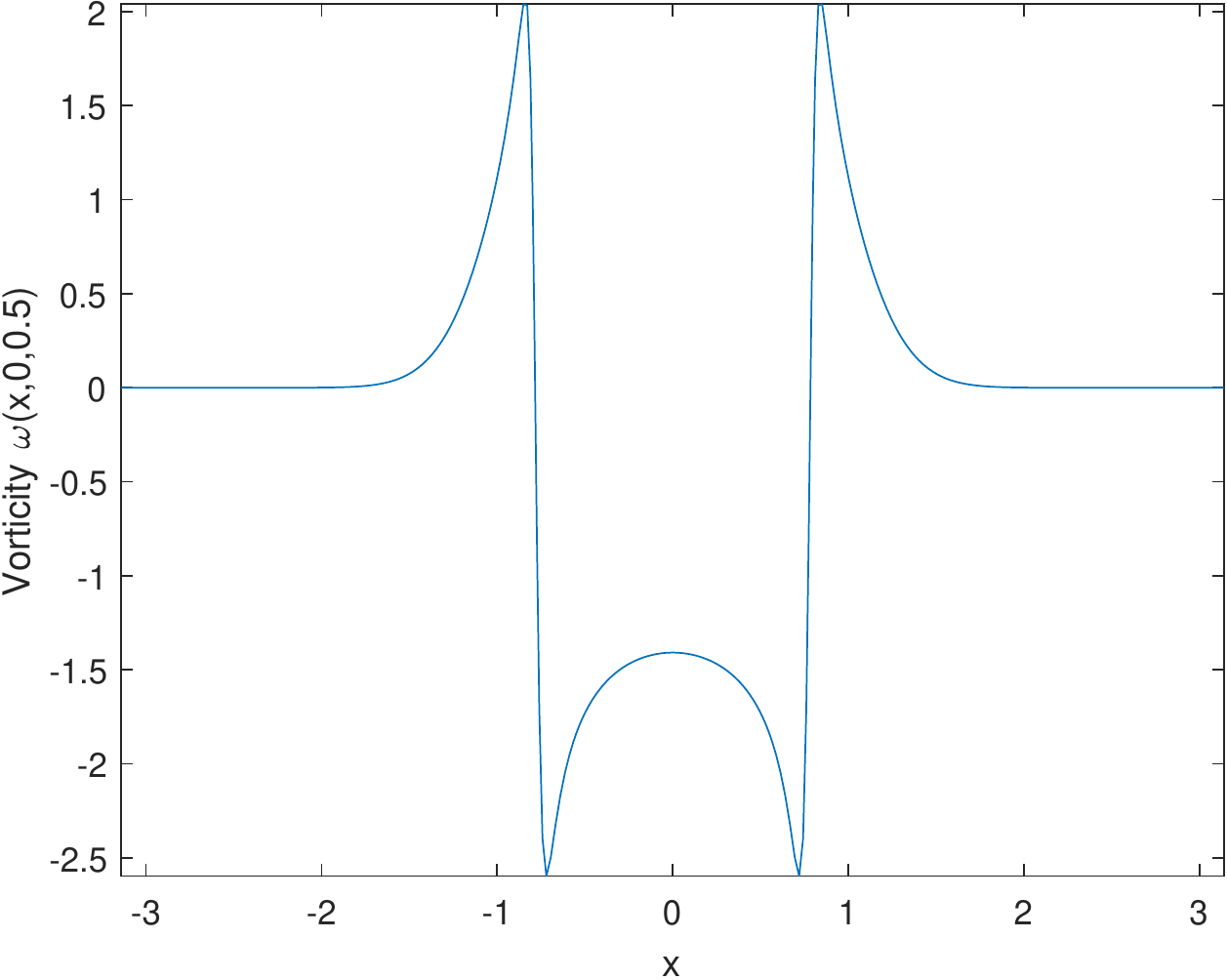}
   \includegraphics[width=.3\textwidth]{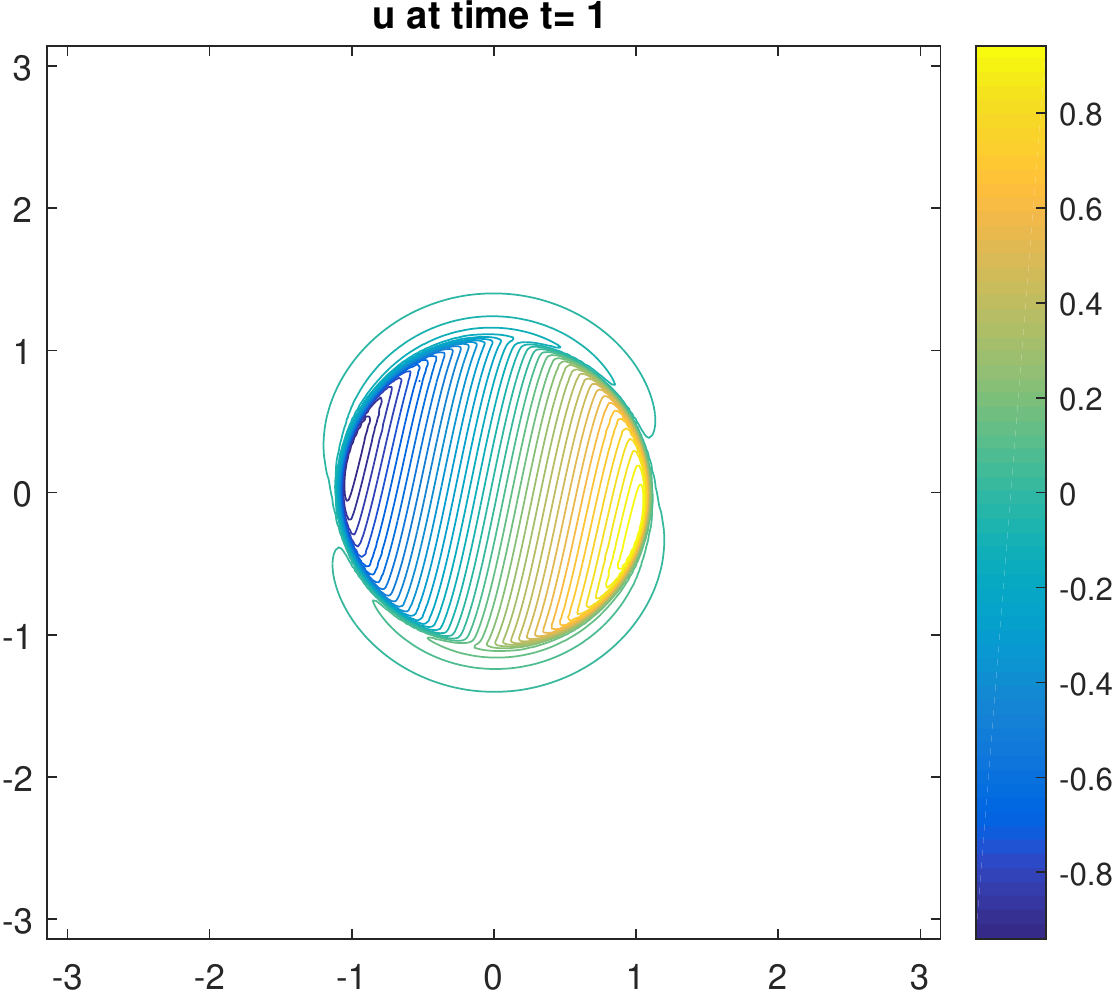}
   \includegraphics[width=.3\textwidth]{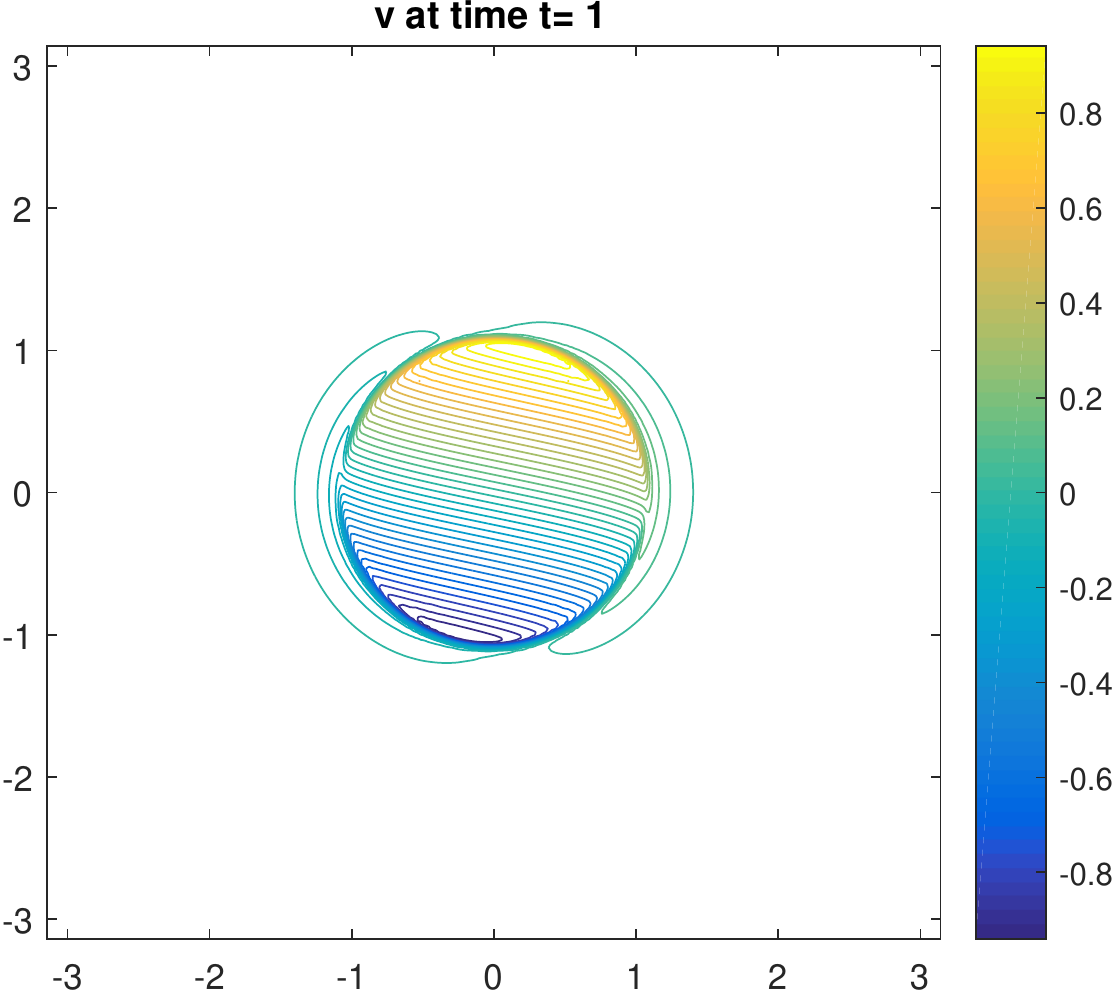}
   \includegraphics[width=.335\textwidth]{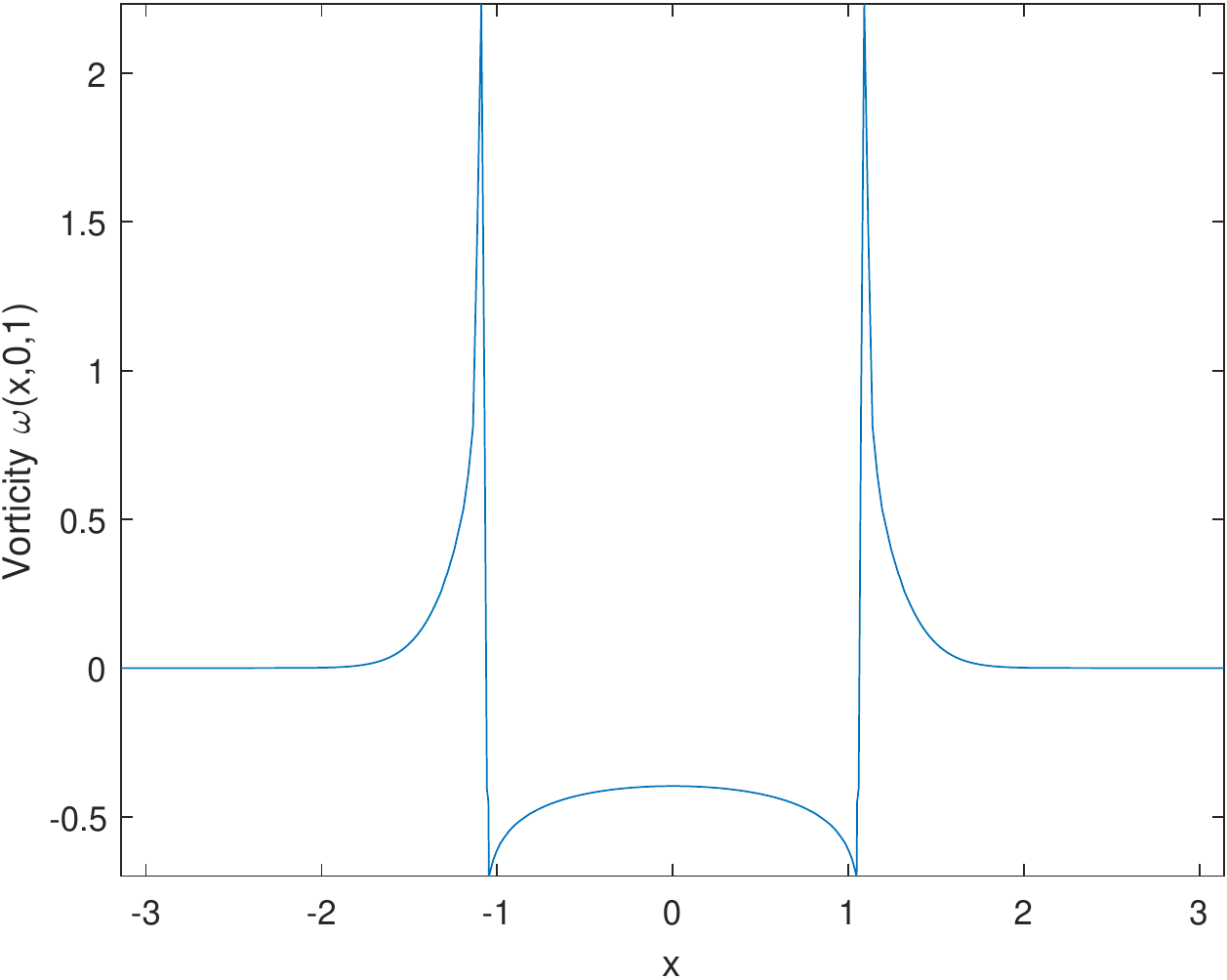}
   \caption{The four plots on the left show the velocities in the $x$ and $y$ directions at times $t = .5$ and $t=1$.  We can see the fluid rotating and expanding.  The two plots in the third column show the vorticity at these times.  We see sharp gradients near the edge of the rotating fluid.\label{fig:Burger1}}
\end{center}
\end{figure}
The first problem we solve uses the initial condition $\mathbf{u} = 5 [-y,x]^T \exp(-3r^2)$ and the boundary conditions are taken to be no-slip boundary conditions on all sides. We solve the problem using $24 \times 24$ leafs, $p = 24$, $\varepsilon = 0.005$, and the fifth order ARK method found in \cite{Carpenter_Kennedy_ARK}.  We use a time step of $k = 1/80$ and solve until time $t_{\rm max} = 5$.  The low viscosity combined with the initial condition produces a rotating flow resembling a vortex that steepens up over time.

\begin{figure}[!ht]
 \begin{center}
   \includegraphics[width=.345\textwidth]{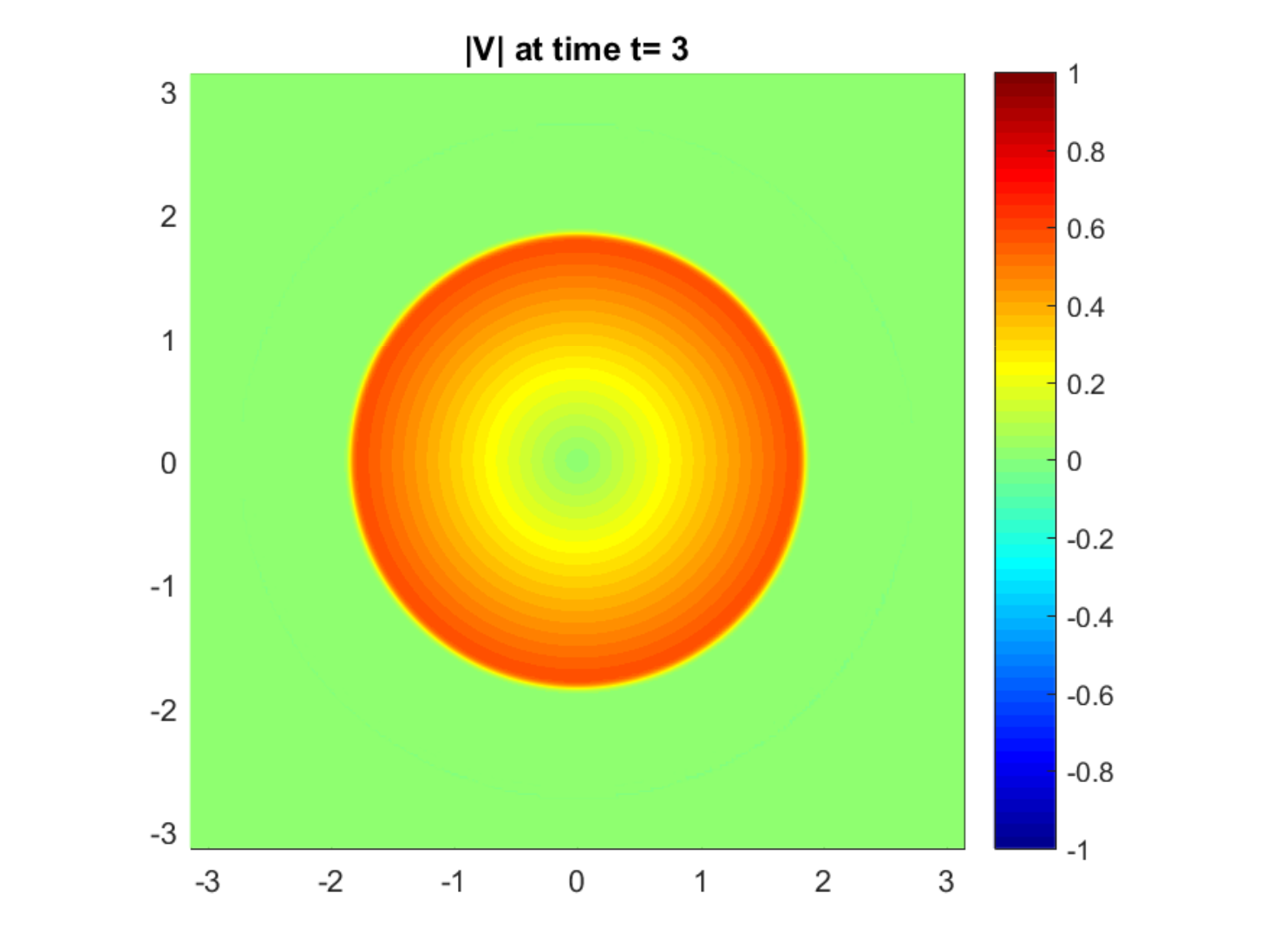}
   \includegraphics[width=.28\textwidth]{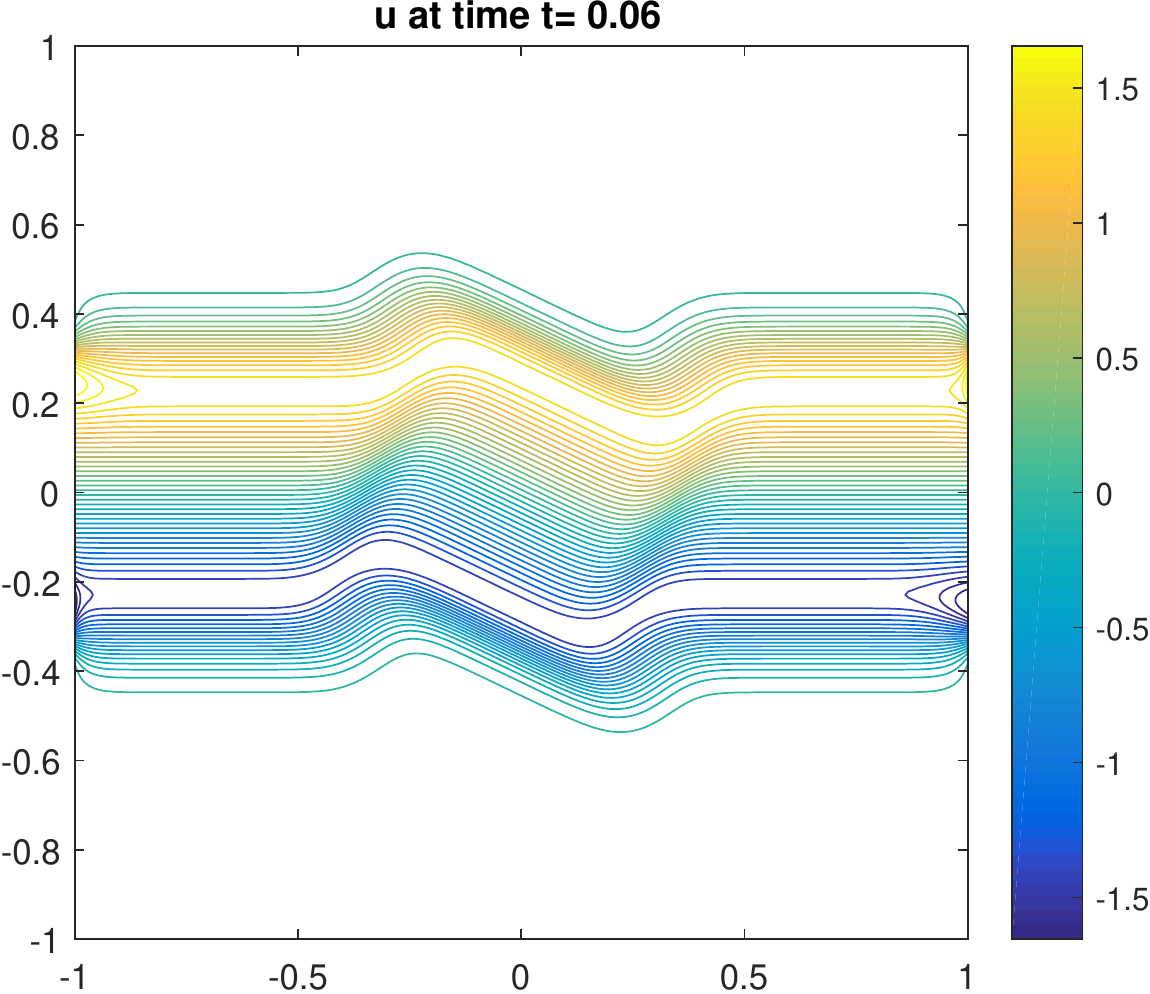}
   \includegraphics[width=.345\textwidth]{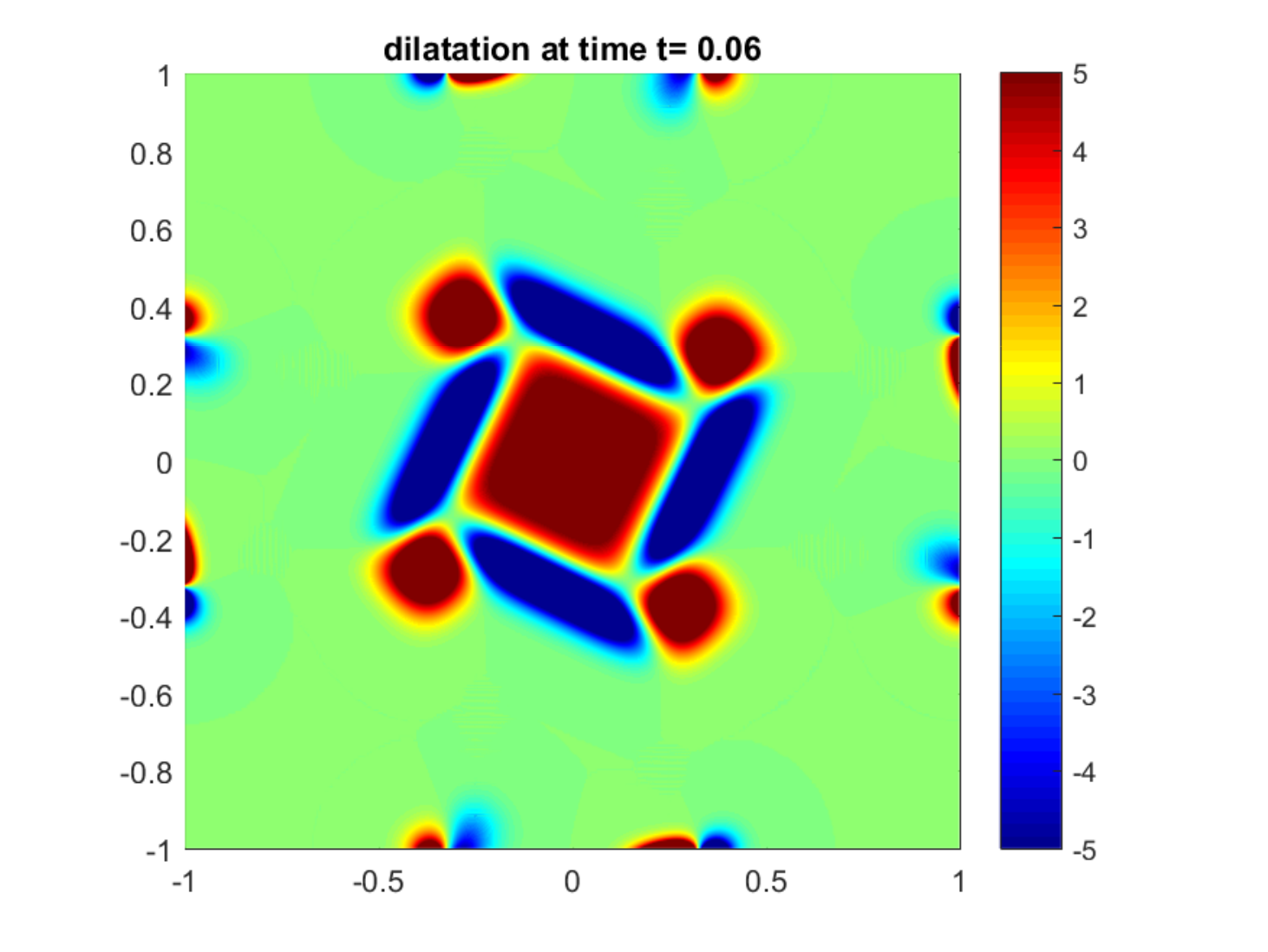}
   \includegraphics[width=.345\textwidth]{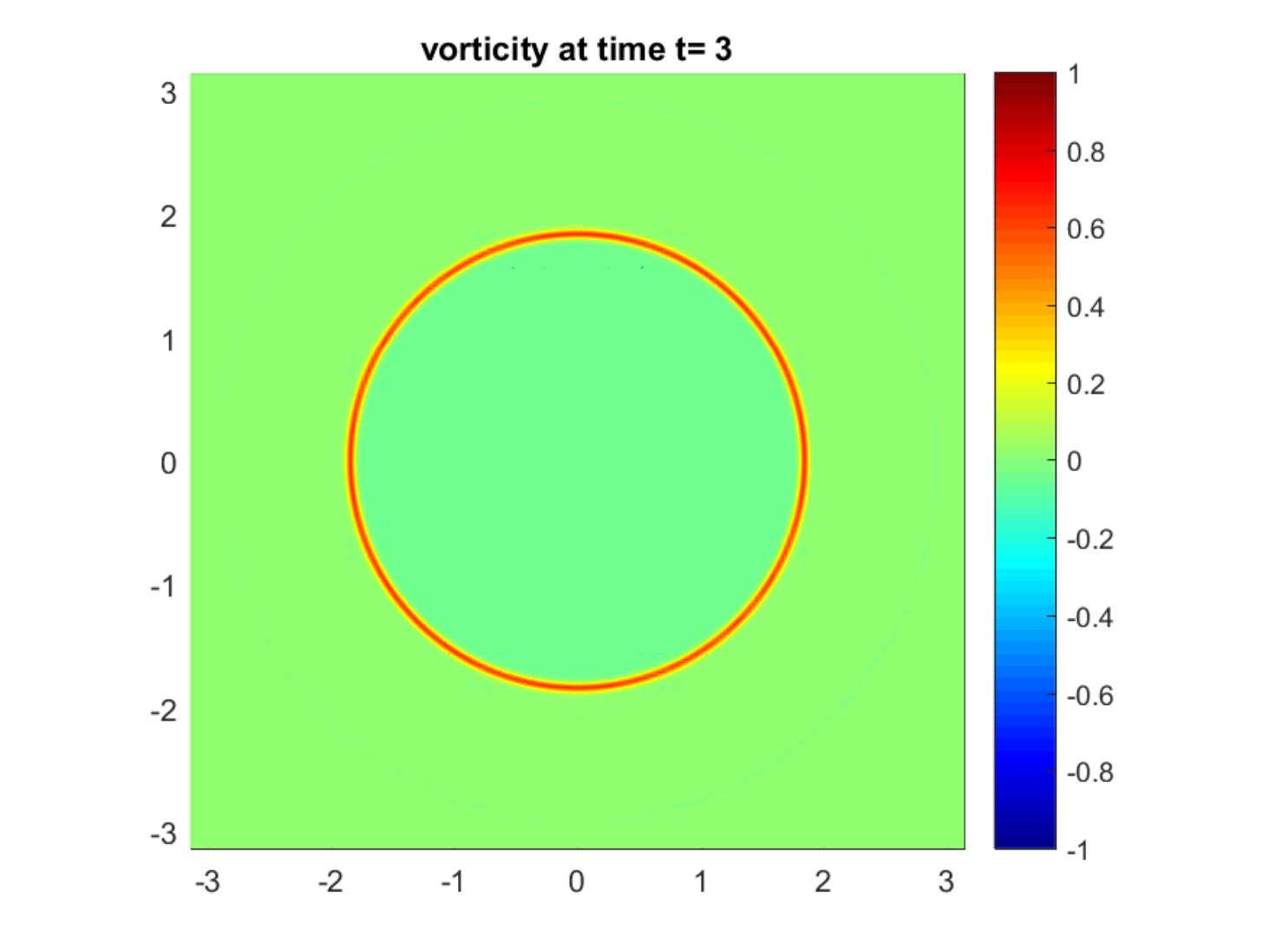}
   \includegraphics[width=.28\textwidth]{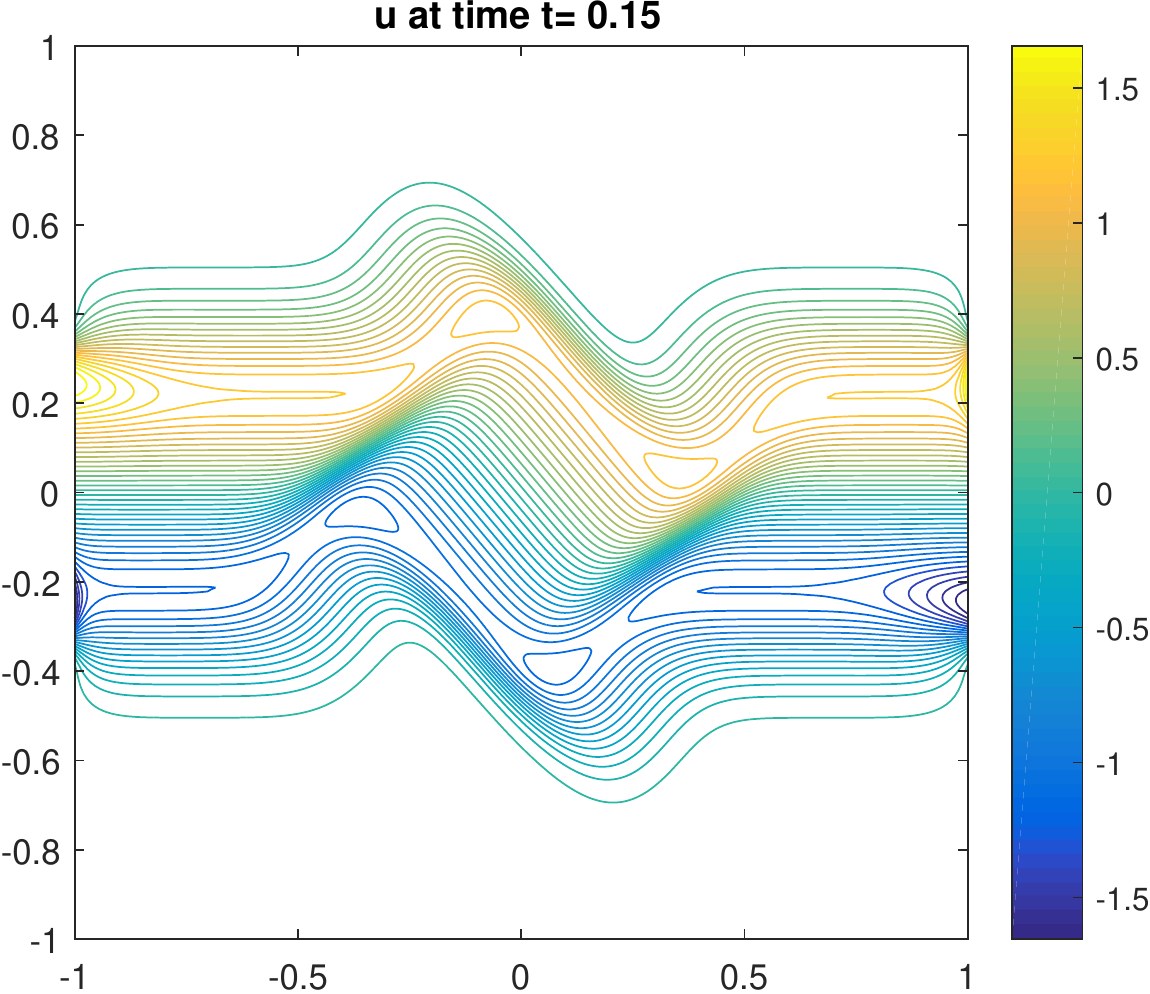}
   \includegraphics[width=.345\textwidth]{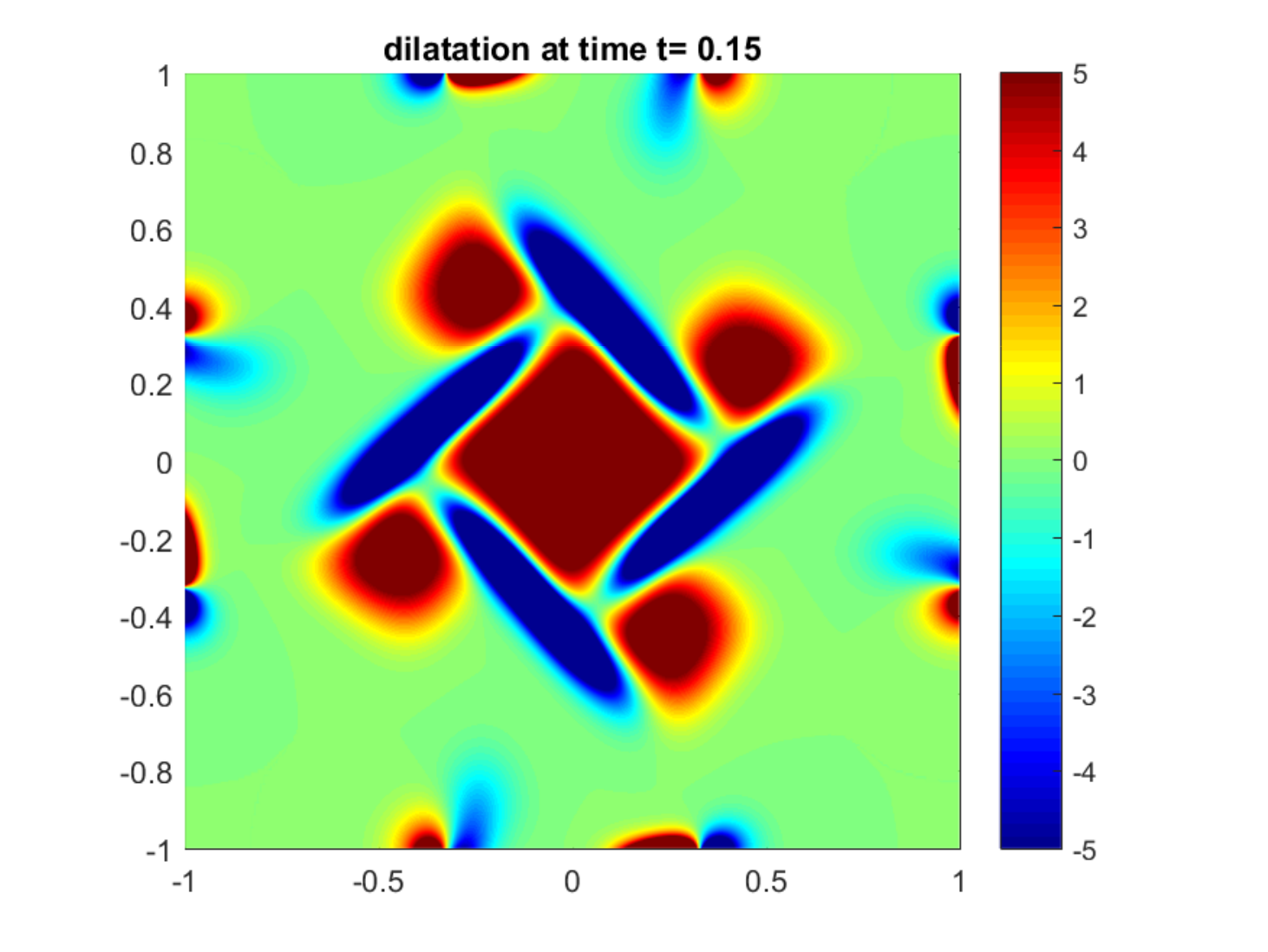}
   \caption{The plots in the first column correspond to the rotating shock problem.  The second and third columns show the velocity in the $x$ direction and the dilatation at times $t = 0.06$ and $t = 0.15$ for the cross flow problem.\label{fig:Burger2}}
\end{center}
\end{figure}

In Figure \ref{fig:Burger1} we can see the velocities at times $t = 0.5$ and $t = 1$. The fluid rotates and expands out and eventually forms a shock like transition. This creates a sharp flow region with large gradients resulting in a flow that may be difficult to resolve with a low order accurate method. These sharp gradients can be seen in the two vorticity plots in Figure \ref{fig:Burger1} along with the speed and vorticity plots in Figure \ref{fig:Burger2}.

In our second experiment we consider a cross stream of orthogonal flows.  We use an initial condition of 
\begin{equation}
\mathbf{u} =  [ 8y\, e^{-36\left(\frac{y}{2}\right)^8},-8x e^{-36\left(\frac{x}{2}\right)^8}]^T,
\end{equation}
and time independent boundary conditions that are compatible with the initial data. 

This initial horizontal velocity drops to zero quickly as we approach $|y| = 0.5$.  For $|y| < 0.5$ the exponential term approaches $\exp(0)$ and the velocity behaves like $u = 8y$.  The flow has changed slightly by $t = 0.06$, but we can see in Figure \ref{fig:Burger2} the flow is moving to the right for $ y > 0$ and the flow is moving the left for $ y < 0$ and all significant behavior is in $|y| < 0.5$.  A plot of the velocity $v$ would show similar behavior.  We also use $24 \times 24$ leafs, $p=24$, $\epsilon = 0.025$, $k = 1/200$, and $t_{\rm max} = 0.75$.  We show plots of the horizontal velocity $u$ and the dilatation at time $t = 0.06$ and $t = 0.15$.  We only show plots before time $t = 0.15$ when the fluid is hardest to resolve and we observe that after $t = 0.15$ the cross streams begin to dissipate.  This problem contains sharp interfaces inside $\mathbf{x} \in [-0.5,0.5]^2$.

\section{Conclusion} \label{sec:conc}
In this two part series we have demonstrated that the spectrally accurate Hierarchial Poincar\'{e}-Steklov solver can be easily extended to handle time dependent PDE problems with a parabolic principal part by using ESDIRK methods. We have outlined the advantages of the two possible ways to formulate implicit Runge-Kutta methods within the HPS scheme and demonstrated the capabilities on both linear and non-linear examples. 

There are many avenues for future work, for example:
\begin{itemize}
\item Extension of the solvers to compressible and incompressible flows.
\item Application of the current solvers to inverse and optimal design problems, in particular for problems where changes in parameters do not require new factorizations. 
\end{itemize}

\bibliography{main_bib_2}
\bibliographystyle{siam}
\end{document}